\documentclass[11pt]{article}
\usepackage{rotating}
\usepackage{amsfonts}
\usepackage{amssymb,amsmath}
\numberwithin{equation}{section}
\usepackage{graphics}
\newcommand{\R}{{\mathbb R}}
\newcommand{\C}{{\mathbb C}}
\newcommand{\N}{{\mathbb N}}

\renewcommand{\d}{\,{\rm d}}            
\newcommand{\D}{{\rm d}}                
\renewcommand{\Re}{\mathop{\mathrm{Re}}}
\renewcommand{\Im}{\mathop{\mathrm{Im}}}
\newcommand{\supp}{\mathop{\mathrm{supp}}}
\newcommand{\dist}{\mathop{\mathrm{dist}}}
\newcommand{\vc}{\mathop{\mathrm{cv}}}
\newcommand{\Rsp}{\mathop{R_{\mathrm{sp}}}}

\newcommand{\DD}{{\cal D}}

\newcommand{\LL}{{\cal L}}

\newcommand{\OO}{{\cal O}}

\newcommand{\RR}{{\cal R}}

\newtheorem{theorem}{Theorem}[section]
\newtheorem{lemma}[theorem]{Lemma}
\newtheorem{definition}[theorem]{Definition}
\newtheorem{proposition}[theorem]{Proposition}
\newtheorem{corollary}[theorem]{Corollary}
\newtheorem{remark}[theorem]{Remark}

\newtheorem{conjecture}[theorem]{Conjecture}
\newtheorem{hypothesis}[theorem]{Hypothesis}

\newcommand{\proof}{{\noindent \bf Proof:\ }}

\def\build#1_#2^#3{\mathrel{
  \mathop{\kern 0pt#1}\limits_{#2}^{#3}}}
\def\eqdef{\buildrel\hbox{\scriptsize{def}}\over =}
\def\QED{\mbox{}\hfill$\Box$}
\def\eps{\epsilon}

\newdimen\texpscorrection
\newdimen\figcenter
\def\figurewithtex #1 #2 #3 #4 #5\cr{\null
  {\goodbreak\figcenter=\hsize\relax
  \advance\figcenter by -#4truecm
  \divide\figcenter by 2
  \begin{figure}[hbt]
  \vskip #3truecm\noindent\hskip\figcenter
  \includegraphics{#1}{\hskip\texpscorrection\input #2 }
  \vskip 0.8truecm{\baselineskip=0.8\baselineskip
  \noindent \vbox{\noindent {\footnotesize #5}}\par}
  \end{figure}}}
\def\point#1 #2 #3 {\rlap{\kern #1 truecm
\raise #2 truecm \hbox{#3}}}

\begin{document}

\title{Spectral asymptotics for large skew-symmetric perturbations
of the harmonic oscillator}

\author{
Isabelle Gallagher \\
Institut de Math\'ematiques de Jussieu\\
Universit\'e de Paris 7\\
Case 7012, 2 place Jussieu\\
75251 Paris Cedex 05, France\\
{\tt gallagher@math.jussieu.fr}\\[3mm]
\and
Thierry Gallay \\ 
Institut Fourier \\
Universit\'e de Grenoble I \\
BP 74\\
38402 Saint-Martin-d'H\`eres, France \\
{\tt Thierry.Gallay@ujf-grenoble.fr}
\and
Francis Nier \\    
IRMAR \\
Universit\'e de Rennes 1 \\
Campus de Beaulieu \\          
35042 Rennes, France \\
{\tt Francis.Nier@univ-rennes1.fr}}
\date{}

\maketitle
\begin{abstract}
Originally motivated by a stability problem in Fluid Mechanics, we study the 
spectral and pseudospectral properties of the differential operator 
$H_\epsilon = -\partial_x^2 + x^2 + i\epsilon^{-1}f(x)$ on $L^2(\R)$, 
where $f$ is a real-valued function and $\epsilon > 0$ a small
parameter. We define $\Sigma(\epsilon)$ as the infimum of the real 
part of the spectrum of $H_\epsilon$, and $\Psi(\epsilon)^{-1}$ as
the supremum of the norm of the resolvent of $H_\epsilon$ along  
the imaginary axis. Under appropriate conditions on $f$, we show 
that both quantities $\Sigma(\epsilon)$, $\Psi(\epsilon)$ go to
infinity as $\epsilon \to 0$, and we give precise estimates of the 
growth rate of $\Psi(\epsilon)$. We also provide an example where
$\Sigma(\epsilon) \gg \Psi(\epsilon)$ if $\epsilon$ is small.
Our main results are established using variational ``hypocoercive''
methods, localization techniques and semiclassical subelliptic 
estimates. 
\end{abstract}

\section{Introduction} \label{intro}

In many evolution equations arising in Mathematical Physics, one
encounters the situation where the generator of the evolution can be
written as a sum of a dissipative and a conservative operator which 
do not commute with each other. In such a case the conservative term 
can affect and sometimes enhance the dissipative effects or the 
regularizing properties of the whole system. For instance, if the 
system has a globally attracting equilibrium, the rate of convergence
towards this steady state can strongly depend on the nature and the
size of the conservative terms. Typical examples illustrating such an 
interplay between diffusion and transport are the kinetic Fokker-Planck 
equation \cite{HerNi}, and the Boltzmann equation \cite{DV}; see also 
\cite{Vi} for a comprehensive study of these phenomena at a more 
abstract level.

In this paper we study a simple linear system which fits into this 
general framework. Given a small parameter $\epsilon > 0$ and a 
smooth, bounded function $f : \R \to \R$, we consider the differential 
operator
\begin{equation}\label{Heps}
  H_\epsilon \,=\, -\partial_x^2 + x^2 + \frac{i}{\epsilon}f(x)~,  
  \quad x \in \R~,
\end{equation}
acting on the Hilbert space $X = L^2(\R)$, with domain $\DD = 
\{u \in H^2(\R)\,;\, x^2 u \in L^2(\R)\}$. Clearly~$H_\epsilon$ is a
bounded, skew-symmetric perturbation of the harmonic oscillator 
$H_\infty = -\partial_x^2 + x^2$. Our goal is to compute the decay
rate in time of the solutions to the evolution equation
\begin{equation}\label{Heq}
  \frac{\D u}{\D t} \,=\, -H_\epsilon u~, \quad u(0) = u_0 \in X~.
\end{equation}
As we shall see, the solutions to \eqref{Heq} decay to zero at least
like $e^{-t}$ as $t \to +\infty$, but the actual convergence rate 
strongly depends on the value of $\epsilon$ and the detailed
properties of $f$, due to the interaction between the symmetric
(dissipative) and the skew-symmetric (conservative) part of the 
generator $-H_\epsilon$. 

Our initial motivation for this study is a specific problem in Fluid
Mechanics which we now briefly describe. As is explained in
\cite{GW1,GW2}, to investigate the long-time behavior of solutions
to the two-dimensional Navier-Stokes equation, it is convenient to
use the the vorticity formulation. In self-similar variables, the 
system reads:
\begin{equation}\label{2Dvort}
  \frac{\partial \omega}{\partial t} + u \cdot \nabla\omega \,=\,
  \Delta \omega + \frac12 x\cdot\nabla \omega + \omega~, \quad
  x \in \R^2~, \quad t \ge 0~,
\end{equation}
where $\omega(x,t) \in \R$ is the vorticity distribution and $u(x,t)
\in \R^2$ is the divergence-free velocity field obtained from $\omega$
via the Biot-Savart law. Equation \eqref{2Dvort} has a family of
stationary solutions, called {\it Oseen vortices}, of the form $\omega
= \alpha G$ where $G(x) = (4\pi)^{-1}e^{-|x|^2/4}$ and $\alpha \in \R$
is a free parameter (the circulation Reynolds number). It turns out
that the linearization of \eqref{2Dvort} at $\alpha G$ has the same
form as \eqref{Heq}, in the sense that the generator can be written as
a difference $\LL - \alpha\Lambda$, where~$\LL$ is a self-adjoint
operator in the weighted space $L^2(\R^2,G^{-1}\d x)$ and $\Lambda$ is
a skew-symmetric perturbation. The analogy goes even further if we
conjugate $\LL - \alpha\Lambda$ with the Gaussian weight~$G^{1/2}$ and
if we neglect a nonlocal, lower-order term in the perturbation
$\Lambda$. The linearized operator then becomes
\begin{equation}\label{Hvort}
  \widetilde H_\alpha  \,=\, -\Delta + \frac{|x|^2}{16} - \frac12 +
  \alpha \tilde f(x)\partial_\theta~, \quad x \in \R^2~,
\end{equation}
where $\partial_\theta = x_1 \partial_2 - x_2\partial_1$ and
$\tilde f(x) = (2\pi |x|^2)^{-1}(1 - e^{-|x|^2/4})$. The operator 
$H_\epsilon$ in \eqref{Heps} is a one-dimensional analog of 
$\widetilde H_\alpha$, and the limit $\epsilon \to 0$ corresponds to
the fast rotation limit $\alpha \to \infty$. Remark that, in this
particular example, the function $\tilde f$ has a unique critical point 
located at the origin, and decreases to zero like $|x|^{-2}$ as 
$|x| \to \infty$. 

The aim of this paper is to study the spectral and pseudospectral
properties of the linear operator $H_\epsilon$ in the limit $\epsilon
\to 0$. Besides the specific motivations explained above, this
question has its own interest from a mathematical point of view, and
turns out to be relatively complex. We have to deal with a
non-self-adjoint problem of (almost) semiclassical type which exhibits
a competition between various microlocal models at different scales,
depending on the structure of the critical points and the decay rate
at infinity of the function $f$. In particular, unlike in the
self-adjoint case, or even in some non-self-adjoint problems such as the
kinetic Fokker-Planck operator (see \cite{HelNi,HerNi,HSS}), the
pseudospectral estimates are not monotone with respect to the
imaginary part of the spectral parameter. Nevertheless, the model
\eqref{Heps} is simple enough so that the analysis can be pushed quite
far, and we believe that our results give a good idea of the phenomena
that can be expected to occur in more general situations. We also
mention that the spectral theory of non-self-adjoint operators, especially
in the semiclassical limit, is a topic of current interest
\cite{Tr,Da,DSZ,Sjo}.

We start with a few basic observations concerning the operator
$H_\epsilon$. As is well known the limiting operator $H_\infty =
-\partial_x^2 + x^2$ is self-adjoint in $L^2(\R)$ with compact
resolvent, and its spectrum is a sequence of simple eigenvalues
$\{\lambda_n^0 \}_{n \in \N}$, where $\lambda_n^0 = 2n+1$. By
classical perturbation theory~\cite{Ka}, it follows that $H_\epsilon$
has a compact resolvent for any $\epsilon > 0$, and that its spectrum is
again a sequence of (simple) eigenvalues $\{\lambda_n(\epsilon)\}_{n
  \in \N}$, with $\Re(\lambda_n(\epsilon)) \to +\infty$ as $n \to
\infty$.  Moreover, it is clear from \eqref{Heps} that the {\it
  numerical range}
\begin{equation}\label{Thetadef}
  \Theta(H_\epsilon) \,=\, \Bigl\{\langle H_\epsilon u,u\rangle_{L^2} 
  \in \C\,;\, u \in \DD\,,~ \|u\|_{L^2} = 1 \Bigr\}
\end{equation}
is contained in the region $\RR_\epsilon \subset \C$ defined by
\begin{equation}\label{RRdef}
  \RR_\epsilon \,=\, \Bigl\{\lambda \in \C \,;\, \Re(\lambda) \ge
  1\,,~\epsilon \Im(\lambda) \in \overline{f(\R)}\Bigr\}~.
\end{equation}
In particular we have $\lambda_n(\epsilon) \in \RR_\epsilon$ for all 
$n \in \N$ and all $\epsilon > 0$. The operator $H_\epsilon$ is 
therefore sectorial, and $H_\epsilon - 1$ is maximal accretive, 
see Section~V-3.10 in \cite{Ka}. By the Lumer-Phillips theorem, 
it follows that
\[
  \|e^{-tH_\epsilon}\| \,\le\, e^{-t}~,\quad \hbox{for all }
  t \ge 0~.
\]
Here and in the sequel we denote by $\|A\|$ the norm of a bounded 
linear operator $A$ in $L^2(\R)$, and by $\Rsp(A)$ its spectral
radius.

Since we want to compute the decay rate of the solutions to 
\eqref{Heq}, we are interested in locating the real part of the 
spectrum of $H_\epsilon$. We thus define:
\begin{equation}\label{Sigmadef}
  \Sigma(\epsilon) \,=\, \inf \Re(\sigma(H_\epsilon)) \,=\, 
  \min_{n \in \N} \Re(\lambda_n(\epsilon))~.
\end{equation}
We also introduce a related pseudospectral quantity, which will play 
an important role in this work:
\begin{equation}\label{Psidef}
  \Psi(\epsilon) \,=\, \left(\sup_{\lambda \in \R} \|(H_\epsilon 
  -i\lambda)^{-1}\|\right)^{-1}~.
\end{equation}
It is easy to see that $\Sigma(\epsilon) \ge \Psi(\epsilon) \ge 1$. 
Indeed, for any $\lambda \in \R$, we have
\begin{equation}\label{specbounds}
  \frac{1}{\dist(i\lambda,\sigma(H_\epsilon))} \,=\, 
  \Rsp((H_\epsilon - i\lambda)^{-1}) \,\le\, 
  \|(H_\epsilon - i\lambda)^{-1}\| \,\le\, 
  \frac{1}{\dist(i\lambda,\Theta(H_\epsilon))}~,
\end{equation}
see e.g. Problem~III-6.16 and Theorem~V-3.2 in \cite{Ka}. Taking the
supremum over $\lambda \in \R$ and using the inclusion $\Theta(H_\epsilon)
\subset \RR_\epsilon$, we thus find
\[
  \frac{1}{\Sigma(\epsilon)} \,\le\, \frac{1}{\Psi(\epsilon)} \,\le\, 
  \sup_{\lambda \in \R}\,\frac{1}{\dist(i\lambda,\RR_\epsilon)} 
 \,=\, 1~,
\]
which is the desired result. Elementary relations between
$\Sigma(\epsilon)$, $\Psi(\epsilon)$ and the norm of the semigroup
$e^{-tH_\epsilon}$ can also be deduced from the following general lemma:

\begin{lemma}\label{elem}
Let $A$ be a maximal accretive operator in a Hilbert space $X$, with
numerical range contained in the sector $\{z\in \C\,;\,|\arg z| \le
\frac{\pi}{2}-2\alpha\}$ for some $\alpha \in (0,\frac{\pi}{4}]$. 
Assume that $A$ is invertible and let 
$$
  \Sigma \,=\, \inf\Re(\sigma(A)) \,>\, 0~, \quad \hbox{and} \quad
  \Psi \,=\, \left(\sup_{\lambda\in\R} \|(A-i\lambda)^{-1}\|
  \right)^{-1}~.
$$
Then the following holds:\\[1mm]
\noindent\textbf{i)} If there exist $C \ge 1$ and $\mu > 0$ such that 
$\|e^{-tA}\| \le C\,e^{-\mu t}$ for all $t \ge 0$, then
\[
  \Sigma \,\ge\, \mu~, \quad \hbox{and} \quad 
  \Psi \,\ge\, \frac{\mu}{1+\log(C)}~\cdotp
\]
\noindent\textbf{ii)} For any $\mu \in (0,\Sigma)$, we have
$\|e^{-tA}\| \le C(A,\mu) \,e^{-\mu t}$ for all $t \ge 0$, where
$$
  C(A,\mu) \,=\, \frac{1}{\pi\tan\alpha}\Bigl(\mu N(A,\mu) + 
  2\pi\Bigr)~, \quad \hbox{and}\quad N(A,\mu) \,=\, 
  \sup_{\lambda \in \R}\|(A-\mu-i\lambda)^{-1}\|~.
$$
\noindent\textbf{iii)} 
If moreover $\mu \in (0,\Psi)$, the quantity $N(A,\mu)$ is not larger 
than $(\Psi-\mu)^{-1}$.
\end{lemma}

For completeness, we give a proof of this lemma in
Appendix~\ref{appendix}. We also mention that a similar result holds
for some maximal accretive operators with $\alpha=0$, in which case
the constant $C(A,\mu)$ has a different expression, see \cite{HerNi}
or \cite{HelNi} for a detailed study of Fokker-Planck operators which
enter in this category. In the case when $A = H_\epsilon$, we have of
course $\Sigma = \Sigma(\epsilon)$, $\Psi = \Psi(\epsilon)$, and the
angle $\alpha$ satisfies $\tan(2\alpha) = \epsilon
\|f\|_{L^{\infty}}^{-1}$, see \eqref{Thetadef}, \eqref{RRdef}. Thus
$\alpha = \OO(\epsilon)$ as $\epsilon \to 0$.

Lemma~\ref{elem} can be used in particular to illustrate the
pseudospectral nature of the quantity $\Psi(\epsilon)$ defined in
\eqref{Psidef}. For semiclassical operators of the form $P_\epsilon =
p(x,\epsilon \partial_x)$, the pseudospectrum is usually defined as 
the set of all $z \in \C$ such that $\epsilon^{N} \|(P_{\epsilon}-z)^{-1}\|
\to \infty$ as $\epsilon \to 0$ for all $N \in \N$, see e.g.
\cite{Pra}. Instead of defining the pseudospectrum as a subset of
$\C$, we find it convenient to introduce here a more flexible notion: 

\begin{definition}
\label{de.pseudo}  
Let $(\omega_{\epsilon})_{\epsilon\in (0,1]}$ be a family of
complex domains, i.e. $\omega_{\epsilon}\subset \C$ for all
$\epsilon\in (0,1]$.\\
We say that $\omega_{\epsilon}$ {\em meets the pseudospectrum} of 
$H_{\epsilon}$ as $\epsilon \to 0$ if 
\[
  \lim_{\epsilon\to 0}\epsilon^{N}\sup_{z\in \omega_{\epsilon}}
  \|(H_{\epsilon}-z)^{-1}\| \,=\, +\infty~, \quad \hbox{for all }
  N\in \N~.
\]
On the contrary, we say that $\omega_{\epsilon}$ {\em avoids the 
pseudospectrum} of $H_{\epsilon}$ as $\epsilon \to 0$  if there
exists $N \in \N$ such that
\[
  \sup_{z\in \omega_{\epsilon}} \|(H_{\epsilon}-z)^{-1}\| \,=\, 
  \OO(\epsilon^{-N})~, \quad \hbox{as } \epsilon \to 0~.
\] 
\end{definition}

If we apply this definition in the case when $\omega_\epsilon = 
\{z \in \C\,;\,\Re(z) \le \mu_{\epsilon}\}$, we arrive at very 
different conclusions depending on whether $\mu_\epsilon \ll 
\Psi(\epsilon)$ or $\mu_\epsilon \gg \Psi(\epsilon)$. Indeed, 
using Lemma~\ref{elem} and Definition~\ref{de.pseudo}, we obtain 
the following result, whose proof is again postponed to 
Appendix~\ref{appendix}: 

\begin{lemma}\label{pseudo}~\\[1mm]
\textbf{i)} For any $\kappa\in (0,1)$, the domain $\{\Re(z) \le 
\kappa\Psi(\epsilon)\}$ avoids the pseudospectrum of $H_{\epsilon}$ 
as $\epsilon\to 0$.\\[1mm]
\textbf{ii)} If $\mu_\epsilon \gg \Psi(\epsilon)(1+\log\Psi(\epsilon) + 
  \log(\epsilon^{-1}))$ in the sense that the ratio goes to $+\infty$
as $\epsilon\to 0$, then the domain $\{\Re(z) \le \mu_\epsilon\}$
meets the pseudospectrum of $H_{\epsilon}$ as $\epsilon\to 0$.
\end{lemma}

\medskip The main purpose of the present work is to investigate the
behavior of both quantities~$\Sigma(\epsilon)$ and~$\Psi(\epsilon)$
as~$\epsilon \to 0$.  In particular, we shall give precise estimates
on $\Psi(\epsilon)$ under rather general assumptions on the function
$f$. Computing $\Sigma(\epsilon)$ is a more delicate task, and we
shall restrict ourselves to a specific example which shows that
$\Sigma(\epsilon)$ can be much larger than $\Psi(\epsilon)$ when
$\epsilon$ is small. This can be expected because $H_\epsilon$ is
highly non-self-adjoint in this regime (we recall that~$\Sigma = \Psi$
for self-adjoint operators).

If the function $f : \R \to \R$ is not identically constant, it is
rather straightforward to verify that $\Sigma(\epsilon) \ge
\Psi(\epsilon) > 1$ for all $\epsilon > 0$, see
Lemma~\ref{nonconstant} below. In view of Lemma~\ref{elem}, we
conclude that the decay properties of the semigroup $e^{-tH_\epsilon}$
are always {\em enhanced} by the skew-symmetric perturbation
$i\epsilon^{-1}f(x)$, if $f$ is not constant. It is therefore 
very natural to ask under which conditions on $f$ the decay rate 
can be made arbitrarily large by choosing $\epsilon$ sufficiently 
small. If we use the pseudospectral quantity $\Psi(\epsilon)$ to 
measure the decay, the complete answer is given by the following 
proposition, which is proved in Section~\ref{compact}. 

\begin{proposition}\label{weakprop}
Assume that $f \in L^\infty(\R) \cap C^0(\R)$. If all level sets of 
$f$ have empty interior, then $\Psi(\epsilon) \to \infty$ as $\epsilon 
\to 0$. In the converse case, $\Psi(\epsilon)$ is uniformly bounded 
for all $\epsilon > 0$.
\end{proposition}
We suspect that Proposition~\ref{weakprop} remains true if we replace
$\Psi(\epsilon)$ by $\Sigma(\epsilon)$, but at the present time there
is no proof that $\Sigma(\epsilon)$ remains bounded if $f$ is constant
on some nonempty open interval.

We now turn to more quantitative versions of
Proposition~\ref{weakprop}, which specify the growth rate of
$\Psi(\epsilon)$ or $\Sigma(\epsilon)$ as $\epsilon \to 0$.  Our first
result in this direction will be obtained using a variational method
recently developed by C.~Villani under the name of ``hypocoercivity''
\cite{Vi}. This approach applies (in particular) to linear operators
of the form $L = A^*A + B$ in a Hilbert space $X$, where $B$ is {\em
skew-symmetric}. Under certain conditions, it allows to compare the
spectral properties of $L$ with those of the self-adjoint operator
$\hat L = A^*A + C^*C$, where $C = [A,B]$.  It is clear that our
problem fits into this general framework: if we set $X = L^2(\R,\C)$,
$A = \partial_x + x$, and~$B = (i/\epsilon)f(x)$, we see that
$H_\epsilon - 1 = A^*A + B$. Since $C = [A,B] = (i/\epsilon)f'(x)$,
the associated self-adjoint operator $\hat H_\epsilon$ defined by
$\hat H_\epsilon - 1 = A^*A + C^*C$ has the explicit form
\begin{equation}\label{hatHeps}
  \hat H_\epsilon \,=\, -\partial_x^2 + x^2 + 
  \frac{1}{\epsilon^2}f'(x)^2~, \quad x \in \R~.
\end{equation}
The point is that the spectral properties of $\hat H_\epsilon$ in the
limit $\epsilon \to 0$ are rather easy to study using semiclassical
techniques, see e.g. Lemma~\ref{hatH} below. Once the spectrum of
$\hat H_\epsilon$ is known, Villani's method allows to deduce useful
information on the original operator $H_\epsilon$. A great interest of
this approach is that it can be easily adapted to nonlinear problems
\cite{Vi,Vi2}.

Following closely the proof of the ``basic hypocoercivity result'' 
in \cite[Section~4]{Vi}, we obtain our first main result:

\begin{theorem}\label{thm1}
Assume that $f \in C^3(\R)$ satisfies $f'',f''' \in L^\infty(\R)$, 
and that there exist $M_1 > 0$ and $\nu \in (0,1/2]$ such that
\begin{equation}\label{hatHbdd}
  \langle \hat H_\epsilon u\,,u\rangle_{L^2} \,=\,  
  \int_\R \Bigl(|\partial_x u|^2 + x^2 |u|^2 + \frac{1}{\epsilon^2}
  \,f'(x)^2 |u|^2\Bigr)\d x \,\ge\, \frac{M_1}{\epsilon^{2\nu}}
  \,\|u\|_{L^2}^2~,
\end{equation}
for all $u \in \DD$ and all $\epsilon \in (0,1]$. Then there exists 
$M_2 > 0$ such that, for all $\epsilon \in (0,1]$, 
\begin{equation}\label{SigmaPsibound}
  \Sigma(\epsilon) \,\ge\, \frac{M_2}{\epsilon^\nu}~, \quad 
  \hbox{and} \quad \Psi(\epsilon) \,\ge\, \frac{M_2}{\epsilon^\nu 
  \log(2/\epsilon)}~\cdotp
\end{equation}
\end{theorem}

Remark that we do not suppose in Theorem~\ref{thm1} that $f$ is
bounded, but only that $f'' $ and~$f'''$ are. Under these more
general assumptions, the perturbation $H_\epsilon-H_\infty$ is
relatively bounded with respect to $H_\infty$, but not necessarily
relatively compact. To test the conclusion on a concrete example,
consider the simple case when $f(x) = x^2$. Then we have~$H_\epsilon =
-\partial_x^2 + (1+i/\epsilon)x^2$, so that~$\Sigma(\epsilon) =
\Re((1+i/\epsilon)^{1/2}) = \OO(\epsilon^{-1/2})$ as $\epsilon \to 0$.
On the other hand, $\hat H_\epsilon = -\partial_x^2 +
(1+4/\epsilon^2)x^2$, hence $\inf(\sigma(\hat H_\epsilon)) =
(1+4/\epsilon^2)^{1/2}$ so that \eqref{hatHbdd} is satisfied with $\nu
= 1/2$. Thus the lower bound~\eqref{SigmaPsibound} for
$\Sigma(\epsilon)$ is optimal, and one can also check that the
estimate for $\Psi(\epsilon)$ is optimal except for the logarithmic
correction. However, as we shall see, the results given by 
Theorem~\ref{thm1} are not always optimal, especially when $f'(x)$
decays algebraically as $|x| \to \infty$.

Motivated by our original example~\eqref{Hvort}, we next consider
a specific class of nonlinearities with a prescribed behavior 
at infinity. 

\begin{hypothesis}\label{hypf}
We assume that $f \in C^3(\R,\R)$ has the following properties:\\[1mm]
i) All critical points of $f$ are non-degenerate; i.e., $f'(x) = 0$  
implies $f''(x) \neq 0$.\\[1mm]
ii) There exist positive constants $C$ and $k$ such that, for all 
$x \in \R$ with $|x| \ge 1$,  
\begin{equation}\label{fdecay}
  \left|\partial_x^\ell\Bigl(f(x) - \frac{1}{|x|^k}\Bigr)\right|  
  \,\le\, \frac{C}{|x|^{k+\ell+1}}~, \quad \hbox{for } \ell = 0,1,2,3~. 
\end{equation}
\end{hypothesis}

Loosely speaking, we consider Morse functions which are bounded
together with their derivatives up to third order, and which behave
like $|x|^{-k}$ as $|x| \to \infty$. This is not the most general
class that we can treat with our methods, but it already contains many
interesting examples, including the original problem \eqref{Hvort}.
The assumption that $f'(x)$ decays algebraically as $|x| \to \infty$
is very important: as we shall see, the decay rate of the semigroup
$e^{-tH_\epsilon}$ is determined not only by the structure of the
critical points of $f$, but also by the behavior of $f$ near infinity,
and under hypothesis \eqref{fdecay} there is a competition between
both regimes which is one of the main interests of our problem.

Under Hypothesis~\ref{hypf}, it is rather straightforward to estimate
the lowest eigenvalue of the self-adjoint operator $\hat H_\epsilon$.
The result is:

\begin{lemma}\label{hatH}
If $f$ satisfies Hypothesis~\ref{hypf}, there exists $M_3 \ge 1$ such 
that, for all $\epsilon \in (0,1]$, 
\begin{equation}\label{hatHbound}
  \frac{1}{M_3\,\epsilon^{2\nu}} \,\le\, \inf \sigma(\hat H_\epsilon) 
  \,\le\, \frac{M_3}{\epsilon^{2\nu}}~, \quad \hbox{where} \quad 
  \nu \,=\, \frac{1}{k+2} ~\cdotp
\end{equation}
\end{lemma}

For the reader's convenience, we sketch the proof of this result in
Appendix~\ref{appendix}. Since \eqref{hatHbound} implies
\eqref{hatHbdd} with the same value of $\nu$, Theorem~\ref{thm1} shows
that $\Sigma(\epsilon), \Psi(\epsilon)$ obey the lower bounds
\eqref{SigmaPsibound} with $\nu = (k+2)^{-1}$ if $f$ satisfies
Hypothesis~\ref{hypf}. It turns out that these estimates are not
optimal in the present case, neither for $\Sigma(\epsilon)$ nor for
$\Psi(\epsilon)$. As is discussed in Section~\ref{variational} below,
this is not an intrinsic limitation of the variational method, but the
proof of a general result such as Theorem~\ref{thm1} fails to take
into account all specificities of $f$. Our second main result gives
optimal estimates on $\Psi(\epsilon)$ under Hypothesis~\ref{hypf}:

\begin{theorem}\label{thm2}
If $f$ satisfies Hypothesis~\ref{hypf}, there exists $M_4 \ge 1$ such 
that, for all $\epsilon \in (0,1]$, 
\begin{equation}\label{Psioptimal}
  \frac{1}{M_4\,\epsilon^{\bar \nu}} \,\le\, \Psi(\epsilon) \,\le\, 
  \frac{M_4}{\epsilon^{\bar \nu}}~, \quad \hbox{where} \quad \bar \nu 
  \,=\, \frac{2}{k+4}~\cdotp
  \end{equation}
\end{theorem}

Remark that $\bar \nu > \nu$, so that the lower bound in
\eqref{Psioptimal} is strictly better than in \eqref{SigmaPsibound}.
In this paper, we shall give {\em two different proofs} of this lower
bound. The first one consists in adapting the variational method that
we already used in the proof of Theorem~\ref{thm1}. This approach
gives the lower bound in~\eqref{Psioptimal} with a logarithmic
correction, as in \eqref{SigmaPsibound}. The second proof, which
removes the logarithmic correction and also provides the upper bound
in \eqref{Psioptimal}, uses rather classical techniques in the
analysis of partial differential operators. The idea is first to
reduce the problem to a bounded domain using a dyadic partition of
unity in $\R$, and then to apply semiclassical subelliptic techniques
as developed, for instance, in \cite{DSZ} and in Chap.~27 of
\cite{Hor}.

\medskip Finally, we discuss the behavior of the spectral quantity
$\Sigma(\epsilon)$. Under Hypothesis~\ref{hypf}, Theorem~\ref{thm2}
already provides a good lower bound on $\Sigma(\epsilon)$, via the
inequality $\Sigma(\epsilon) \ge \Psi(\epsilon)$. In fact, if one
believes that the analysis of small random perturbations of
non-self-adjoint semiclassical operators presented in \cite{HS} can be
carried over to the operator $H_{\epsilon}$, one also
expects that the lower bound $\Sigma(\epsilon) \ge \Psi(\epsilon)$ is
optimal for a generic $f$ in a wide class of functions. However
$\Sigma(\epsilon)$ can be much larger than $\Psi(\epsilon)$ in some
particular situations. This is the case for instance when $f(x) = x$,
because one can show that $\Sigma(\epsilon) = \OO(\epsilon^{-2})$ and
$\Psi(\epsilon) = \OO(\epsilon^{-2/3})$, see Section~\ref{fx} below.
To conclude this paper, we consider a more interesting example which
satisfies Hypothesis~\ref{hypf} and is inspired by the original
problem \eqref{Hvort}.

\begin{proposition}\label{pr.speclb}
Fix $k > 0$ and assume that 
\begin{equation}\label{fex}
  f(x) \,=\, \frac{1}{(1+x^2)^{k/2}}~, \quad x \in \R~.
\end{equation}
Then there exists a constant $M_5 > 0$ such that the lowest real part
of the spectrum satisfies, for all $\epsilon \in (0,1]$,  
\begin{equation}\label{Sigmaconj}
  \Sigma(\epsilon) \,\ge\, \frac{M_5}{\epsilon^{\nu'}}~, \quad
  \hbox{where}\quad \nu' \,=\, \min\Bigl\{\frac12\,,\,
  \frac{2}{k+2}\Bigr\}~.
\end{equation}
\end{proposition}

Since $\bar \nu < \nu'$, it follows from \eqref{Psioptimal}
and \eqref{Sigmaconj} that $\Sigma(\epsilon)/\Psi(\epsilon) \to
\infty$ as $\epsilon \to 0$. In particular, if we take 
$\mu_\epsilon = \frac12 M_5\,\epsilon^{-\nu'}$, Lemma~\ref{pseudo} implies 
that the domain $\{\Re(z) \le \mu_\epsilon\}$ meets the pseudospectrum
of $H_\epsilon$ as $\epsilon \to 0$ in the sense of 
Definition~\ref{de.pseudo}. It is also interesting to test the 
conclusions of Lemma~\ref{elem} in this situation: we know from 
\textbf{ii)} that $\|e^{-tH_\epsilon}\| \le C\,e^{-\mu_\epsilon t}$ 
for some $C \ge 1$, but it follows from \textbf{i)} that $C\ge 
e^{-1}\exp(\mu_\epsilon/\Psi(\epsilon))$. Thus, although we can prove that
the semigroup $e^{-tH_\epsilon}$ decays with an exponential rate 
$\mu_\epsilon \gg \Psi(\epsilon)$, this estimate is inadequate
for small times because the prefactor $C$ is exponentially large. 
In contrast, if we choose $\mu_\epsilon \le \kappa\Psi(\epsilon)$
for some $\kappa < 1$, we can take $C = \OO(\epsilon^{-1})$ by 
Lemma~\ref{elem}-\textbf{iii)}.

Proposition~\ref{pr.speclb} will be proved in Section~\ref{spectral}
below, using a complex deformation method (which exploits the
analyticity properties of $f$) and the same localization techniques 
as in the proof of Theorem~\ref{thm2}. We shall also give heuristical 
arguments, supported by accurate numerical computations, which show 
that the lower bound in \eqref{Sigmaconj} is {\em optimal} in the 
sense that the exponent $\nu'$ cannot be improved. 

\medskip The rest of this paper is organized as follows. In
Section~\ref{compact} we prove some general properties of
$\Psi(\epsilon)$ using mainly compactness arguments.
Section~\ref{variational} is devoted to the proof of
Theorem~\ref{thm1} using the approach developed by Villani \cite{Vi}.
We also show how that the variational method can be adapted to give
the lower bound on $\Psi(\epsilon)$ in Theorem~\ref{thm2}, up to a
logarithmic correction.  In Section~\ref{resolvent}, we give accurate
bounds on the resolvent $(H_\epsilon - i\lambda)^{-1}$ for $\lambda
\in \R$ using a dyadic decomposition of the real axis and
semiclassical subelliptic estimates. In particular we prove
Theorem~\ref{thm2}, and we also provide a family
$(\omega_{\epsilon})_{\epsilon\in (0,1]}$ of complex domains which
avoid the pseudospectrum of $H_{\epsilon}$ as $\epsilon\to 0$
according to Definition~\ref{de.pseudo}, with some non trivial
$\epsilon$-dependent geometry. Finally, in Section~\ref{spectral}, the
lower bound \eqref{Sigmaconj} on $\Sigma(\epsilon)$ is proved and
illustrated by numerical computations. The proofs of
Lemmas~\ref{elem}, \ref{pseudo}, and \ref{hatH} are collected in
Appendix~\ref{appendix}.

\medskip\noindent{\bf Acknowledgements.} The authors would
like to thank V. Bach, Y. Colin de Verdi\`ere, M. Hitrik, A. Joye, 
G. Perelman, C. Villani, S. Vu Ngoc, and M. Zworski for fruitful 
discussions.


\section{Compactness estimates} 
\label{compact}

In this section we establish a few general properties of the 
quantities $\Sigma(\epsilon)$ and $\Psi(\epsilon)$ defined 
in \eqref{Sigmadef}, \eqref{Psidef}. We first show that 
any nontrivial function $f$ does affect the spectrum of the 
operator $H_\epsilon$.

\begin{lemma}\label{nonconstant}
Assume that $f \in L^\infty(\R)$ is not a
constant. Then $\Sigma(\epsilon) \ge \Psi(\epsilon) > 1$ for 
all $\epsilon > 0$. 
\end{lemma}

\proof Suppose on the contrary that $\Psi(\epsilon) = 1$ for
some $\epsilon > 0$. Then, in view of \eqref{Psidef}, 
\eqref{specbounds}, there exists $\lambda_0 \in \R$ such that
$\|(H_\epsilon - i\lambda_0)^{-1}\| = 1$. Thus, for any $n \in 
\N^*$, we can find $\phi_n \in L^2(\R)$ with $\|\phi_n\| = 1$ 
such that
\[
  1 - \frac{1}{n} \,\le\, \|(H_\epsilon - i\lambda_0)^{-1}\phi_n\|
  \,\le\, 1~.
\]
Let $\psi_n = (H_\epsilon - i\lambda_0)^{-1}\phi_n$, so that $\psi_n
\in \DD$ and $\phi_n = (H_\epsilon - i\lambda_0)\psi_n$. For all 
$n \in \N^*$ we have
\[
  \|\psi_n'\|^2 + \|x\psi_n\|^2 \,=\, \Re \langle (H_\epsilon - i
  \lambda_0)\psi_n,\psi_n \rangle \,=\,\Re \langle \phi_n,\psi_n 
  \rangle \,\le\, 1~.
\]
By Rellich's theorem, the sequence $\{\psi_n\}$ is compact 
in $L^2(\R)$. Thus, after extracting a subsequence, we can assume 
that $\psi_n$ converges in $L^2(\R)$ to some limit $\psi$. 
By construction $\|\psi\| = 1$ and
\[
  \|\psi'\|^2 + \|x\psi\|^2 \,\le\, \liminf_{n \to \infty}
  (\|\psi_n'\|^2 + \|x\psi_n\|^2) \,\le\, 1~. 
\]
Since $\|\psi'\|^2 + \|x\psi\|^2 = \langle H_\infty \psi,\psi\rangle$, 
this inequality implies that $\psi(x) = \pi^{-1/4} \,e^{-x^2/2}$ 
(up to a constant phase factor). Moreover, since $\psi_n \to \psi$
in $L^2(\R)$ and $\|(H_\epsilon - i\lambda_0)\psi_n\| \le 1$ for all
$n$, we necessarily have $\|(H_\epsilon - i\lambda_0)\psi\| \le 1$. 
But $(H_\epsilon - i\lambda_0)\psi = \psi + i(\epsilon^{-1}f-\lambda_0)
\psi$ and therefore
\[
  \|(H_\epsilon - i\lambda_0)\psi\|^2 \,=\, \int_\R |\psi(x)|^2
  \Bigl\{1 + \Bigl(\frac{1}{\epsilon}f(x)-\lambda_0\Bigr)^2 
  \Bigr\}\d x \,=\, 1 + \int_\R |\psi(x)|^2
  \Bigl(\frac{1}{\epsilon}f(x) -\lambda_0\Bigr)^2\d x~,
\]
hence it is possible to have $\|(H_\epsilon - i\lambda_0)\psi\| \le 1$
only if $f(x) = \epsilon\lambda_0$ almost everywhere. \QED

\bigskip
Using similar arguments, we next prove that $\Psi(\epsilon) \to 
\infty$ as $\epsilon \to 0$ if and only if the level sets of 
$f$ have empty interior (we assume here, for simplicity, that
$f$ is continuous).

\bigskip\noindent{\bf Proof of Proposition~\ref{weakprop}:}
Assume first that $\Psi(\epsilon)^{-1}$ does not converge to
zero as $\epsilon \to 0$. Then, according to definition 
\eqref{Psidef}, there exist a positive constant $\delta$, 
a sequence $\{\epsilon_n\}$ of positive numbers, a sequence 
$\{\lambda_n\}$ of real numbers, and a sequence $\{u_n\}$ of
normalized vectors in $L^2(\R)$ such that $\epsilon_n \to 0$
as $n \to \infty$ and 
\[
  \|(H_{\epsilon_n} -i\lambda_n)^{-1} u_n\| \,\ge\, \delta \,>\, 0~,
  \quad \hbox{for all } n \in \N~.
\]
Let $v_n = (H_{\epsilon_n} -i\lambda_n)^{-1} u_n$, so that 
$v_n \in \DD$ and $u_n = (H_{\epsilon_n} -i\lambda_n)v_n$. 
The sequence $\{v_n\}$ is bounded in $L^2(\R)$, because
$\|v_n\| \le \|u_n\| \le 1$ for all $n \in \N$. Moreover, 
we have the identity
\begin{equation}\label{unvn}
  \langle u_n,v_n\rangle \,=\, \|v_n'\|^2 + \|xv_n\|^2 + 
  \frac{i}{\epsilon_n}\,\langle fv_n,v_n\rangle - i\lambda_n
  \|v_n\|^2~.
\end{equation}
Taking the real parts of both sides, we find $\|v_n'\|^2 + \|xv_n\|^2 
\le |\langle u_n,v_n\rangle| \le 1$ for all $n \in \N$. By Rellich's
theorem, the sequence $\{v_n\}$ is therefore compact in $L^2(\R)$. 
Thus, after extracting a subsequence, we can assume that $v_n$ 
converges in $L^2(\R)$ to some limit $v$. By construction, 
$v \in H^1(\R)$ and $\|v\|_{L^2} \ge \delta > 0$. 

On the other hand, multiplying both sides of \eqref{unvn} by 
$\epsilon_n$ and taking imaginary parts, we see that 
\[
  \langle fv_n,v_n\rangle - \epsilon_n \lambda_n \|v_n\|^2 
  \,=\, \epsilon_n \Im \langle u_n,v_n\rangle 
  \,\xrightarrow[n\to \infty]{}\, 0~.
\]
Since $\langle fv_n,v_n\rangle \to \langle fv,v\rangle$, we conclude
that
\[
  \lim_{n \to \infty} \epsilon_n \lambda_n \,=\, 
  \frac{\langle fv,v\rangle}{\|v\|^2} ~\eqdef~ \mu~.  
\]
Finally, for any test function $\phi \in C_0^\infty(\R)$, we
have
\begin{align*}
  \epsilon_n \langle u_n,\phi\rangle \,&=\, \epsilon_n 
  \langle (H_{\epsilon_n} -i\lambda_n)v_n,\phi\rangle \\
  \,&=\, \epsilon_n \langle v_n,H_\infty \phi \rangle + i\langle 
  fv_n,\phi \rangle -i \epsilon_n \lambda_n \langle v_n,\phi\rangle~.
\end{align*}
Taking the limit $n \to \infty$, we obtain $\langle fv,\phi\rangle
= \mu \langle v,\phi\rangle$ for all $\phi \in C_0^\infty(\R)$, 
which implies that $(f-\mu)v=0$. Since $f$ and $v$ are continuous 
functions and $v \not\equiv 0$, this is possible only if $f^{-1}(\mu)$ 
has non-empty interior. Thus we have shown that, if all level sets of
$f$ have empty interior, then necessarily $\Psi(\epsilon) \to 
\infty$ as $\epsilon \to 0$. 

Conversely, assume that there exists $\mu \in \R$ such that 
$f^{-1}(\mu)$ contains a nonempty open interval $I$. Choose 
$\phi \in C_0^\infty(I)$ such that $\phi \not\equiv 0$, and let
\[
  \psi \,=\, H_\infty \phi \,\equiv\, (H_\epsilon -\frac{i\mu}{\epsilon})
  \phi \qquad (\hbox{for any }\epsilon > 0)~.
\]
Observe that $\psi \not\equiv 0$, because $0 \notin
\sigma(H_\infty)$. Thus
\[
  \Psi(\epsilon)^{-1} \,\ge\, \|(H_\epsilon -\frac{i\mu}{\epsilon})^{-1}\|
  \,\ge\, \frac{\|\phi\|}{\|\psi\|}~,
\]
hence $\Psi(\epsilon) \le \|\psi\|/\|\phi\|$ for all $\epsilon > 0$. 
\QED
 

\section{Variational estimates} \label{variational} 

The purpose of this section is to adapt the general method developed
in \cite{Vi} to our concrete problem, and to show that it allows to
obtain precise estimates on the semigroup and the resolvent of the
operator $H_\epsilon$. To demonstrate the efficiency of this approach,
we first give a short proof of Theorem~\ref{thm1} in
Section~\ref{proofthm1}. We next consider in more detail three model
problems, and we show how to modify the proof of Theorem~\ref{thm1} so
as to obtain optimal estimates in each case. With this information at
hand, we consider in Section~\ref{optimal} the general situation where
the function $f$ satisfies Hypothesis~\ref{hypf} and we prove the
first half of Theorem~\ref{thm2}, namely the lower bound on
$\Psi(\epsilon)$ up to a logarithmic correction. We mention however
that a complete proof of Theorem~\ref{thm2} will be given in
Section~\ref{resolvent}, using a different method.

\subsection{Proof of Theorem~\ref{thm1}}\label{proofthm1}
This section is devoted to the proof of Theorem~\ref{thm1}. 
Consider any $f \in C^3(\R)$ such that $f''$ and~$f'''$ belong to 
$L^\infty(\R)$, and assume that estimate \eqref{hatHbdd} holds 
for some $\nu \in (0,1/2]$. Let~$u(x,t)$ be a solution to the
parabolic equation \eqref{Heq}, namely
\begin{equation}\label{ueq}
  \partial_t u(x,t) \,=\, \partial_x^2 u(x,t) - x^2 u(x,t)
  -\frac{i}{\epsilon}\,f(x) u(x,t)~.
\end{equation}
To control the evolution of $u(x,t)$, we introduce the quadratic 
functional
\begin{equation}\label{Phidef}
  \Phi(t) \,=\, \int_\R\Bigl(\frac12 |u|^2 + \frac{\alpha}{2}
  (|\partial_x u|^2 + x^2|u|^2) + \beta\Re ((\partial_x \bar u)
  \,i f'(x)u) + \frac{\gamma}{2} f'(x)^2 |u|^2\Bigr)\d x~,
\end{equation}
where $\alpha,\beta,\gamma$ are positive constants to be determined
below. We assume that $4\beta^2 \le \alpha\gamma$, so that
\begin{equation}\label{Philow}
  \Phi(t) \,\le\, \int_\R\Bigl(\frac12 |u|^2 + \frac{3\alpha}{4}
  (|\partial_x u|^2 + x^2|u|^2) + \frac{3\gamma}{4} f'(x)^2 |u|^2
  \Bigr)\d x~.
\end{equation}
Note that, according to the general method introduced in~\cite{Vi}, 
the functional $\Phi$ is a linear combination of the
quantities~$\|u\|^2$, $\|Au\|^2$, $\Re \langle Au,Cu\rangle$
and~$\|Cu\|^2$, where as in the introduction~$A = \partial_x + x$
and~$C = [A,B]$ with~$B = (i/\eps) f(x)$.  To compute the time
derivative of $\Phi(t)$, we use the identities:
\begin{align}\label{id1}
&\frac12 \frac{\D}{\D t}\int_\R |u|^2\d x \,=\, 
  \Re \int_\R \bar u \partial_t u\d x \,=\, -\int_\R 
  (|\partial_x u|^2 + x^2 |u|^2)\d x~, \\ \label{id2}
&\frac12 \frac{\D}{\D t}\int_\R \alpha (|\partial_x u|^2 + x^2|u|^2)
  \d x \,=\, -\int_\R \alpha |\partial_x^2 u - x^2 u|^2\d x 
  -\Re \int_\R \alpha (\partial_x \bar u)\frac{i}{\epsilon}
  f'(x)u \d x~, \\  \nonumber
&\frac{\D}{\D t}\,\Re\int_\R \beta (\partial_x \bar u)\,i f'(x)u \d x
  \,=\, -\frac{1}{\epsilon}\int_\R \beta f'(x)^2 |u|^2 \d x
  -\Re\int_\R \beta\bar u \,i f'''(x) \partial_x u \d x \\  \label{id3}
&\hspace{51mm} +\,2\Re\int_\R \beta (\partial_x \bar u)\,if'(x)
  (\partial_x^2 u - x^2 u)\d x~, \\ \nonumber
&\frac12 \frac{\D}{\D t}\int_\R \gamma  f'(x)^2|u|^2\d x \,=\,
  - \int_\R \gamma f'(x)^2 (|\partial_x u|^2 + x^2 |u|^2)\d x \\
& \hspace{42mm} - 2\Re \int_\R \gamma f'(x)f''(x)\bar u\partial_x
  u\d x \label{id4}, 
\end{align}
which follow easily from \eqref{ueq} after some integrations 
by parts. To estimate the various terms in the expression 
of $\Phi'(t)$, we denote
\[
  K_j \,=\, \sup_{x \in \R}|\partial_x^j f(x)|~, \quad 
  j = 2,3~,
\]
and we use the following bounds: 
\begin{align}
& -\alpha \Re \int_\R (\partial_x \bar u)\frac{i}{\epsilon}
  f'(x)u \d x \,\le\, \frac14 \int_\R |\partial_x u|^2\d x 
  + \frac{\alpha^2}{\epsilon^2} \int_\R f'(x)^2 |u|^2 \d x~, \label{id11}\\
& -\beta \Re \int_\R \bar u \,i f'''(x) \partial_x u \d x
  \,\le\, K_3 \beta \int_\R |u \partial_x u|\d x \,\le\, 
  K_3 \beta \int_\R (|\partial_x u|^2 + x^2 |u|^2)\d x~,\label{id22}\\
& 2\beta\Re\int_\R (\partial_x \bar u)\,if'(x)
  (\partial_x^2 u{-}x^2 u)\d x 
  \,\le\, \frac{\alpha}{2}\int_\R 
  |\partial_x^2 u{-}x^2 u|^2\d x + \frac{2\beta^2}{\alpha}
  \int_\R f'(x)^2 |\partial_x u|^2 \d x ,\label{id33}\\
& - 2\gamma \Re \int_\R f'(x)f''(x)\bar u\partial_x u\d x
  \,\le\, \frac14 \int_\R |\partial_x u|^2\d x + 
  4 \gamma^2 K_2^2 \int_\R f'(x)^2 |u|^2\d x~.\label{id44}
\end{align}
We thus obtain
\begin{align}
  \Phi'(t) \,&\le\, 
  \Bigl(-1 + \frac14 + \frac14 + K_3\beta\Bigr)
  \int_\R (|\partial_x u|^2 + x^2 |u|^2)\d x \nonumber \\
  \,&+\, \Bigl(-\alpha + \frac{\alpha}2\Bigr)\int_\R |\partial_x^2 u - 
  x^2 u|^2\d x + \Bigl(-\frac{\beta}{\epsilon} + 
  \frac{\alpha^2}{\epsilon^2} + 4\gamma^2 K_2^2\Bigr)\int_\R 
  f'(x)^2 |u|^2 \d x \nonumber\\
  \,&+\, \Bigl(-\gamma + \frac{2\beta^2}{\alpha}\Bigr)
  \int_\R f'(x)^2 (|\partial_x u|^2 + x^2 |u|^2)\d x~. \label{estimatephi'}
\end{align}
We now fix
\begin{equation}\label{betadef}
  \beta \,=\, \min\Bigl(\frac{1}{4K_3}\,,\,\frac{1}{32K_2}\Bigr)~, 
  \quad \alpha \,=\, \Bigl(\frac{\beta\epsilon}{4}\Bigr)^{1/2}~,
  \quad \gamma \,=\, 8 \Bigl(\frac{\beta^3} \epsilon\Bigr)^{1/2}~.
\end{equation}
Then $4\beta^2 \le \alpha\gamma$, $4\gamma^2 K_2^2 \le \beta/(4\epsilon)$,  
and we arrive at the simpler estimate
\begin{align}\label{Phidot}
  \Phi'(t) \,\le\, &-\frac14 \int_\R (|\partial_x u|^2 + x^2 |u|^2)\d x
  -\frac{\beta}{2\epsilon} \int_\R f'(x)^2 |u|^2 \d x \\ \label{dummy}
  &- \frac{\alpha}{2}\int_\R |\partial_x^2 u - x^2 u|^2\d x 
  - \frac{\gamma}{2}\int_\R f'(x)^2 (|\partial_x u|^2 + x^2 |u|^2)\d x~.
\end{align}
In what follows we neglect both terms in \eqref{dummy} and 
keep only the upper bound \eqref{Phidot}. Combining \eqref{Phidot}
with \eqref{hatHbdd} and assuming (without loss of generality)
that $\epsilon \le 2\beta$, we arrive at
\[
  \Phi'(t) \,\le\,  -\frac18 \int_\R (|\partial_x u|^2 + x^2 |u|^2)\d x
  -\frac{\beta}{4\epsilon} \int_\R f'(x)^2 |u|^2 \d x 
  - \frac{M_1}{8}\Bigl(\frac{2\beta}{\epsilon}\Bigr)^\nu
  \int_\R |u|^2\d x~.
\]
Using finally \eqref{Philow}, we conclude that $\Phi'(t) \le 
-\eta\,\Phi(t)$, where
\begin{equation}\label{choiceeta}
  \eta \,=\, \min\left(\frac{1}{6\alpha}\,,\,
  \frac{\beta}{3\epsilon\gamma}\,,\,\frac{M_1}{4}
  \Bigl(\frac{2\beta}{\epsilon}\Bigr)^\nu
  \right)~.
\end{equation}
Since $\nu \le 1/2$ and $\alpha,\beta,\gamma$ are given by
\eqref{betadef}, there exists $M_2 > 0$ such that 
$\eta \ge 2M_2 \epsilon^{-\nu}$ for all $\epsilon \in (0,1]$.

What we have done so far is to construct a quadratic functional $\Phi$
which is equivalent (with $\epsilon$-dependent constants) to the
square of the norm in $D(H_\infty^{1/2}) = \{u \in H^1(\R)\,;\, 
xu \in L^2(\R)\}$, and which satisfies $\Phi(t) \le
e^{-\eta (t-s)} \Phi(s)$ for all $t \ge s \ge 0$. To deduce some 
information on the semigroup $e^{-tH_\epsilon}$ in $L^2(\R)$, we use
the following standard argument. Let $u(x,t)$ be a solution
to~\eqref{ueq} with initial data $u_0 \in L^2(\R)$. In view of
\eqref{id1}, we have
\[
  \int_0^{\sqrt{\epsilon}} (\|\partial_x u(t)\|^2 + \|xu(t)\|^2)\d t 
  \,=\, \frac12 (\|u_0\|^2 - \|u(t)\|^2) \,\le\, \frac12 \|u_0\|^2~,
\]
hence there exists $\tau \in (0,\sqrt{\epsilon}]$ such that 
\[
  \|u(\tau)\|^2 \,\le\, \|u_0\|^2~, \quad \|\partial_x u(\tau)\|^2 + 
  \|xu(\tau)\|^2 \,\le\, \frac{\|u_0\|^2}{2\sqrt{\epsilon}}~.
\]
Thus, using the definition \eqref{Phidef} of $\Phi$, the values
\eqref{betadef} for $\alpha$, $\beta$, $\gamma$, and the fact that 
$|f'(x)| \le |f'(0)| + K_2 |x|$ for all $x \in \R$, we see that
\[
  \Phi(\sqrt{\epsilon}) \,\le\, \Phi(\tau) \,\le\, \frac{C}{\epsilon}
  \|u_0\|^2~, \quad \hbox{for some } C > 0 \hbox{ independent of }
  \epsilon \in (0,1]~.
\]
As a consequence, for any $t \ge 0$, we obtain
\[
  \|u(t+\sqrt{\epsilon})\|^2 \,\le\, 2\Phi(t+\sqrt{\epsilon})
  \,\le\, 2\,e^{-\eta t}\Phi(\sqrt{\epsilon}) \,\le\,
  \frac{2C}{\epsilon}\,e^{-\eta t}\|u_0\|^2~,
\]
hence
\[
  \|e^{-(t+\sqrt{\epsilon})H_\epsilon}\| \,\le\, \Bigl(\frac{2C}{\epsilon}
  \Bigr)^{1/2} \,e^{-M_2 t/\epsilon^\nu}~, \quad \hbox{for all }t \ge 0~.
\]
Since we also have $\|e^{-tH_\epsilon}\| \le 1$ for all $t \ge 0$, 
it follows from Lemma~\ref{elem}-\textbf{i)} that $\Sigma(\epsilon) 
\ge M_2/\epsilon^{\nu}$. Moreover, proceeding as in the proof of 
that lemma (see Appendix~\ref{appendix}), we find
\begin{align*}
  \Psi(\epsilon)^{-1} \,&\le\, \int_0^{\infty} \|e^{-tH_\epsilon}\|\d t
  \,\le\, \int_0^{\sqrt{\epsilon}}\d t + \int_0^\infty \min\Bigl\{
  1\,,\, \Bigl(\frac{2C}{\epsilon}\Bigr)^{1/2} \,e^{-M_2 t/\epsilon^\nu}
  \Bigr\}\d t \\
  \,&\le\, \sqrt{\epsilon} + \frac{\epsilon^\nu}{M_2} \Bigl(1 + \log
  \Bigl(\frac{2C}{\epsilon}\Bigr)^{1/2}\Bigr)~.
\end{align*}
Since $\nu \le 1/2$, we deduce that $\Psi(\epsilon)^{-1} \le C' 
\epsilon^\nu \log(2/\epsilon)$. This concludes the proof of 
Theorem~\ref{thm1}. \QED

\subsection{Three model problems}\label{model}

A simple and general result such as Theorem~\ref{thm1} cannot give
optimal estimates for all possible choices of the function $f$.  In
particular, if $f$ satisfies Hypothesis~\ref{hypf}, the lower bounds
\eqref{SigmaPsibound} with $\nu = (k+2)^{-1}$ (see Lemma~\ref{hatH})
are too pessimistic.  In this paragraph we consider three
representative examples, and we show in each case that the proof of
Theorem~\ref{thm1} gives in fact an optimal lower bound on
$\Psi(\epsilon)$ provided we choose in an appropriate way the
parameters~$\alpha$, $\beta$, $\gamma$ entering in the
definition~\eqref{Phidef} of the functional~$\Phi$. These examples are
the building blocks which will allow us to treat in
Section~\ref{optimal} the general case of a function $f$ satisfying
Hypothesis~\ref{hypf}. It is interesting to note that the proof of
Theorem~\ref{thm2} in Section~\ref{resolvent} involves a similar
enumeration of cases, see Proposition~\ref{pr.Hpseudoraff}.

\subsubsection{The case of a nondegenerate critical point}\label{fx2}
We first consider the situation where $f$ has a unique, non
nondegenerate critical point, and $|f'|$ is bounded away from 
zero outside a neighborhood of this point. The simplest example of 
such a function is~$f(x) = x^2$. In that case we can just follow
the proof of Theorem~\ref{thm1} and choose the parameters $\alpha$, 
$\beta$, $\gamma$ as in \eqref{betadef}. Indeed, the lower bound
\eqref{hatHbdd} holds with $\nu = 1/2$, hence it follows from
\eqref{choiceeta} that $\Phi'(t) \le -\eta\,\Phi(t)$ with 
$\eta = \OO(\epsilon^{-1/2})$. This result is clearly optimal, 
because we know that $\Sigma(\epsilon) = \Re((1+i/\epsilon)^{1/2})
=  \OO(\epsilon^{-1/2})$ if $f(x) = x^2$. 
 
\subsubsection{The case when $f$  has no critical point}\label{fx}
We next study the case when $|f'|$ is bounded away from zero. 
The simplest paradigm here is~$f(x) = x$. In this example, 
we have
\[
  H_\epsilon \,=\, -\partial_x^2 +x^2 + \frac i \epsilon x \,=\,
  -\partial_x^2 + \Bigl(x+\frac i{2\epsilon}\Bigr)^2 + 
  \frac{1}{4\epsilon^2}~,
\]
hence $\Sigma(\epsilon) = 1 + (4\epsilon^2)^{-1}$. One can also
estimate the resolvent of $H_\epsilon$, and one finds that
$\Psi(\epsilon) =  \OO(\epsilon^{-2/3})$, see
Proposition~\ref{pr.Hpseudoraff}.
To recover this result using the variational method, it is necessary
to modify slightly the proof of Theorem~\ref{thm1}. Indeed, if we
still choose $\alpha$, $\beta$, $\gamma$ as in \eqref{betadef}, then
\eqref{choiceeta} again implies $\eta = \OO(\epsilon^{-1/2})$, because
the lower bound \eqref{hatHbdd} holds with $\nu = 1$. However, since
$f' \equiv 1$, it is more natural here to set $\gamma = 0$ in
\eqref{Phidef} and to assume that $\beta^2 \le \alpha/4$, so that
$\Phi(t)$ is strictly positive. Then proceeding exactly as above and
using the fact that the terms \eqref{id4}, \eqref{id22}, and
\eqref{id44} vanish identically, we see that \eqref{estimatephi'}
becomes
\begin{align*}
  \Phi'(t) \,&\le\, 
  \Bigl(-1 + \frac14   +  \frac{2\beta^2}{\alpha}\Bigr)
  \int_\R (|\partial_x u|^2 + x^2 |u|^2)\d x \nonumber \\
  \,&+\, \Bigl(-\alpha + \frac{\alpha}2\Bigr)\int_\R |\partial_x^2 u - 
  x^2 u|^2\d x \,+\, \Bigl(-\frac{\beta}{\epsilon} + 
  \frac{\alpha^2}{\epsilon^2}\Bigr)\int_\R f'(x)^2 |u|^2 \d x~.
\end{align*}
We need~$\alpha^2 \le \epsilon\beta/2$, so we can choose for instance
\[
  \alpha \,=\, \Bigl(\frac{\epsilon^2}{16}\Bigr)^\frac13~, \quad 
  \hbox{and} \quad \beta \,=\, \Bigl(\frac{\epsilon}{32}\Bigr)^\frac13~.
\]
Then we easily obtain $\Phi'(t) \le -\eta \, \Phi(t)$ with
\[
  \eta \,=\,  \min\left(\frac{1}{3\alpha}\,,\,
  \frac{2\beta}{3\epsilon }\right) \,=\,  \OO(\eps^{-\frac23})~.
\]
{\bf Remark: } In the more general case when $f'$ is not identically 
constant, but still bounded away from zero, we have to deal with 
the additional term which appears in the left-hand side of
\eqref{id22}. Since $|f'''| \le C|f'|$, this term is easily 
bounded as in \eqref{id11}, with $\alpha/\epsilon$ replaced by 
$\beta$. 

\subsubsection{The case of a critical point at infinity}\label{fxk}

Finally, we consider the situation where $f$ has no critical 
point, but $f'(x)$ decreases algebraically as $x \to \pm \infty$. 
We fix $k > 0$ and focus on the following example: 
\[
  f(x) \,=\, \int_0^x\frac{\d y}{(1+y^2)^{\frac{k+1}2}}~,
  \quad x \in \R~,
\]
which nearly satisfies Hypothesis~\ref{hypf}. In particular, $f'(x) 
=  \OO(|x|^{-k-1})$ as $|x| \to \infty$. The proof of Theorem~\ref{thm2}
given in Section~\ref{resolvent} applies to this example, and shows 
that $\Psi(\epsilon) =  \OO(\epsilon^{-\bar \nu})$ where $\bar \nu =
2/(k+4)$. On the other hand, the lower bound \eqref{hatHbdd} holds
with $\nu = (k+2)^{-1} < \bar \nu$, hence the proof of 
Theorem~\ref{thm1} only gives $\Phi'(t) \le -\eta \, \Phi(t)$ with $\eta = 
 \OO(\epsilon^{-\nu})$. However, we shall show how to choose the 
parameters~$\alpha, \beta$, $\gamma$ in \eqref{Phidef} to obtain 
the optimal result in this case too.

So let us again consider the functional~$\Phi(t)$, supposing as
usual~$\beta^2 \leq \alpha\gamma/4$ to ensure \eqref{Philow}, and let
us go through estimates~\eqref{id11} to~\eqref{id44}. The main
advantage we have here compared to the general situation considered in
Theorem~\ref{thm1} is that $|f''| + |f'''| \le C_k|f'|$, where $C_k$ 
denotes (here and below) a generic constant depending only on~$k$. 
Estimates~\eqref{id11} and \eqref{id33}, which involve only~$f'$, are 
unchanged. Estimate~\eqref{id22} is replaced by
\[
  -\beta \Re \int_\R \bar u \,i f'''(x) \partial_x u \d x
 \,\le\, \frac14  \int_\R |\partial_x u|^2  \d x + C_k \beta^2 
 \int_\R f'(x)^2 |u|^2 \d x~.
\]
Finally, after integrating by parts, we bound the left-hand side
of~\eqref{id44} as follows: 
\[
 - 2\gamma \Re \int_\R f'(x)f''(x)\bar u\partial_x u\d x \,=\,
 \gamma \int_\R (f''(x)^2 + f'(x)f'''(x))|u|^2 \d x \,\le\, 
 C_k \gamma \int_\R f'(x)^2 |u|^2\d x~.
\]
Summarizing, estimate \eqref{estimatephi'} becomes
\begin{align*}
  \Phi'(t) \,\le\, &-\frac12 \int_\R (|\partial_x u|^2 + x^2 |u|^2)\d x 
   \,+\,\Bigl(-\frac{\beta}{\epsilon} + \frac{\alpha^2}{\epsilon^2} 
   + C_k (\beta^2 + \gamma)\Bigr)\int_\R f'(x)^2|u|^2 \d x \\
  &-\, \frac{\alpha}2 \int_\R |\partial_x^2 u - x^2 u|^2\d x \,+\, 
  \Bigl(-\gamma + \frac{2\beta^2}{\alpha}\Bigr) \int_\R   
  f'(x)^2 (|\partial_x u|^2 + x^2 |u|^2)\d x~. 
\end{align*}
We now assume that $\alpha^2 \le  \beta\epsilon /4$ and $C_k(\beta^2 
+ \gamma) \le \beta/(4\epsilon)$. Using in addition the lower
bound~\eqref{hatHbdd}, which holds for the value $\nu = (k+2)^{-1}$ 
in view of Lemma~\ref{hatH}, we obtain
\[
  \Phi'(t) \,\le\, -\frac14 \int_\R (|\partial_x u|^2 + x^2 |u|^2)\d x 
  -\frac{\beta}{4\epsilon}\int_\R f'(x)^2|u|^2 \d x
  -\frac{M_1}{4}\Bigl(\frac{\beta}{\epsilon}\Bigr)^\nu
  \int_\R |u|^2 \d x~.
\]
It follows that~$\Phi'(t) \le -\eta\, \Phi(t)$ with 
\[
  \eta \,=\,  \min \left(\frac1{3\alpha}\,,\,\frac{\beta}{3\epsilon
  \gamma}\,,\,\frac{M_1}{2}\Bigl(\frac{\beta}{\epsilon}\Bigr)^\nu
  \right)~.
\]
If we now choose $\alpha = \frac12 \epsilon^\frac{2}{k+4}$, $\beta = 
\epsilon^{-\frac k{k+4}}$, and $\gamma = 8 \epsilon^{-\frac{2k+2}
{k+4}}$, we obtain $\eta =  \OO(\epsilon^{-\frac{2}{k+4}})$ which is the 
desired result.

\subsection{Improved decay rate under Hypothesis~\ref{hypf}}
\label{optimal}
Let us now consider an arbitrary function~$f$ satisfying
Hypothesis~\ref{hypf}. The analysis of Section~\ref{model}
suggests  considering three different regions of space (near the 
critical points of $f$, near infinity, or away from those regions),
and shows how to choose the parameters~$\alpha$, $\beta$, $\gamma$ 
so that the functional~$\Phi$ defined by \eqref{Phidef} decays
with the optimal rate in each region. Since different choices are
needed in different regions, it is very natural to generalize the
definition \eqref{Phidef} and to allow the parameters $\alpha$, 
$\beta$, $\gamma$ to depend on the space variable $x$. So we consider 
again the functional
\[
  \Phi(t) \,=\, \int_\R\Bigl(\frac12 |u|^2 + \frac{\alpha}{2}
  (|\partial_x u|^2 + x^2|u|^2) + \beta\Re ((\partial_x \bar u)
  \,i f'(x)u) + \frac{\gamma}{2} f'(x)^2 |u|^2\Bigr)\d x~,
\]
where now $\alpha$, $\beta$, $\gamma$ are positive functions of~$x$, 
depending also on $\epsilon$. To choose these functions in the
simplest possible way, we first decide that 
\begin{equation}\label{alphagamrel}
  \alpha \,=\, \left(\frac{\beta\epsilon}{4}\right)^\frac12~, 
  \quad \hbox{and} \quad \gamma \,=\, 8\left(\frac{\beta^3}\epsilon
  \right)^\frac12~.
\end{equation}
Indeed this relation between $\alpha, \gamma$ and $\beta,\epsilon$ was
already used in \eqref{betadef}, so if we impose \eqref{alphagamrel}
our set of parameters will be convenient near the critical points of
$f$, provided $\beta$ is small enough. Moreover, as is easily
verified, the parameters we chose in Section~\ref{fxk} also satisfy
\eqref{alphagamrel}, so this definition will be adapted near infinity
provided $\beta$ is of the order of $\epsilon^{-k/(k+4)}$ in that
region. Finally, in the intermediate region we can choose the same
parameters $\alpha$, $\beta$, $\gamma$ as near the critical points of
$f$, because (as was observed in Section~\ref{fx}) the functional
$\Phi$ will decay with the suboptimal rate $\eta =  \OO(\epsilon^{-1/2})$
which is still better than what we have near infinity. Summarizing,
the choice \eqref{alphagamrel} is appropriate in all regions provided
$\beta$ is small near the critical points of $f$ and $\beta =
 \OO(\epsilon^{-k/(k+4)})$ near infinity. Note that \eqref{alphagamrel}
implies $\alpha\gamma = 4\beta^2$, so that the upper bound
\eqref{Philow} is still valid.
\figurewithtex 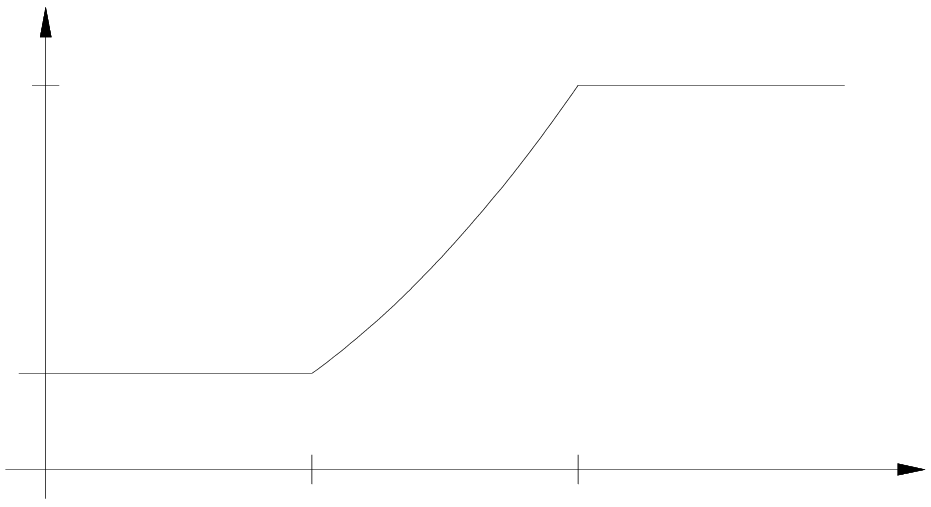 Fig1.tex 5.500 10.500 {\bf Fig.~1:} The
graph of the function $\beta$ defined in \eqref{betadef2}. Here 
$B_\epsilon = A\,\epsilon^{-\frac1{2(k+4)}}$, and $A > 0$ is sufficiently
large so that all critical points of $f$ are contained in the 
interval $[-A+1,A-1]$.\cr

It remains to choose the function $\beta$ appropriately. Let 
$\beta_0 > 0$ be a small constant, and $A > 0$ be large enough 
so that that all critical points of~$f$ are contained in the 
interval $[-A+1,A-1]$. Both parameters are independent of~$\epsilon$, 
and their precise values will be specified later. We define the 
function~$\beta : \R \to \R_+$ by
\begin{equation}\label{betadef2}
  \beta(x) \,=\, \left\{\begin{array}{lll} 
  \beta_0 & \hbox{if} & |x| \le A~, \\[1mm]
  \beta_0{\displaystyle \Bigl(\frac{|x|}{A}\Bigr)^{2k}} & \hbox{if} & 
  |A| \le |x| \le B_\epsilon~, \\[2mm]
  \beta_0\,\epsilon^{-\frac{k}{k+4}}& \hbox{if} & |x| \ge B_\epsilon~,
  \end{array}\right.
\end{equation}
where $k > 0$ is as in Hypothesis~\ref{hypf} and $B_\epsilon =
A\,\epsilon^{-\frac1{2(k+4)}}$. The graph of $\beta$ is represented in
Fig.~1. It is easily verified that $\beta$ has the following useful
properties: for any~$c > 0$, there exists $\epsilon_0 > 0$ such that,
for any $\epsilon \in (0,\epsilon_0)$ and all~$x \in \R$,
\begin{equation}\label{betaprop}
  \epsilon^\frac12 \beta'(x)^2 \,\le\, c \beta(x)^\frac32~, \quad
  \epsilon \beta(x) \,\le\, c~, \quad \hbox{and}\quad  
  \epsilon \beta'(x)^2 \,\le\, c \beta (x)~.  
\end{equation}

We now compute the time derivative of $\Phi(t)$, taking into account
the additional terms involving derivatives of $\alpha$, $\beta$,
$\gamma$ which result from integrations by parts. We do not give
the full details here, but simply indicate the main modifications with 
respect to the corresponding calculations in Section~\ref{proofthm1}. 
Equation \eqref{id1} is of course unchanged, so we start with
\eqref{id2}. Since~$\alpha$ now depends on~$x$, the right-hand side 
of~\eqref{id2} contains an additional term which can be bounded as 
follows:
\[
  -\Re \int_\R \alpha' (\partial_x \bar u)(\partial_x^2 u - x^2 u) \d x 
  \,\le\, \frac14 \int_\R \alpha |\partial_x^2 u - x^2 u|^2 \d x + 
  \int_\R \frac{\alpha'^2}\alpha |\partial_x u|^2 \d x~.
\]
Remark that, in view of \eqref{betaprop}, we have for~$\epsilon > 0$ 
small enough:
\[
  \frac{\alpha'(x)^2}{\alpha(x)} \,=\, \frac{\epsilon^{1/2}\beta'(x)^2}
  {8\beta(x)^{3/2}} \,\le\, \frac{1}{12}~, \quad
  \hbox{for all }x \in \R~. 
\]
On the other hand, since $\alpha^2 = \beta\epsilon/4$, we have by 
\eqref{id11}
\[
  - \Re \int_\R \alpha (\partial_x \bar u)\frac{i}{\epsilon}
  f'(x)u \d x \,\le\, \frac1{4} \int_\R |\partial_x u|^2\d x 
  + \frac{1}{4\epsilon} \int_\R \beta f'(x)^2 |u|^2 \d x~.  
\]
We conclude that
\begin{align}
  \frac12 \frac{\D}{\D t}\int_\R \alpha (|\partial_x u|^2 + x^2|u|^2)
  \d x \,\le\, &- \int_\R \frac{3\alpha}4 |\partial_x^2 u - x^2 u|^2 
  \d x \nonumber\\ \label{id2'}
  &+\frac13 \int_\R |\partial_x u|^2 \d x + \int_\R \frac{\beta}
  {4\epsilon} f'(x)^2 |u|^2 \d x~.
\end{align}

We next consider \eqref{id3}. In that case one should add to the
right-hand side the terms
\[
  \Re \int_\R \beta' \bar u if'(x) (\partial_x^2 u - x^2 u) \d x 
  \,-\, \Re\int_\R \beta' \bar u i f''(x) \partial_x u\d x~.
\]
The first one is easily controlled, for~$\epsilon$ small enough, by
\begin{align*}
  \Re \int_\R \beta' \bar u if'(x)(\partial_x^2 u - x^2 u) \d x  
  \,&\le\, \int_\R \frac\alpha{8} |\partial_x^2 u - x^2 u|^2 \d x 
  + \int_\R \frac{2\beta'^2}\alpha f'(x)^2 |u|^2\d x    \\
  \,&\le\, \int_\R \frac\alpha{8} |\partial_x^2 u - x^2 u|^2 \d x 
  + \int_\R \frac{\beta}{12\epsilon} f'(x)^2 |u|^2 \d x~,
\end{align*}
because $\beta'^2/\alpha \ll \beta/\epsilon$ due to \eqref{alphagamrel}, 
\eqref{betaprop}. Similarly, the second term is estimated by
\begin{align*}
  -\Re \int_\R \beta'\bar u i f''(x) \partial_x u\d x  
  \,&\le\, \frac1{12} \int_\R |\partial_x u|^2 \d x + \int_\R 3 
  \beta'^2 f''(x)^2 |u|^2 \d x \\
  \,&\le\, \frac1{12} \int_\R |\partial_x u|^2 \d x + \int_\R
  \frac\beta{12\epsilon} f'(x)^2 |u|^2 \d x~,
\end{align*}
where in the last inequality we have used~\eqref{betaprop}, as well as
the fact that $|f''(x)| \le C_k |f'(x)|$ for all $x$ in the support
of~$\beta'$. On the other hand, we have as in \eqref{id33}: 
\[
  2\Re \int_\R \beta (\partial_x \bar u)i f'(x)(\partial_x^2 u - x^2 u)
  \d x \,\le\, \int_\R \frac\alpha{2} |\partial_x^2 u - x^2 u|^2 \d x 
  + \int \frac{\gamma}{2} f'(x)^2 |\partial_x u|^2 \d x~,
\]
where we have used the fact that $4\beta^2 = \alpha\gamma$. Finally, 
the term involving $f'''(x)$ in \eqref{id2} is estimated as follows. 
We first note that
\[
  -\Re \int_\R \beta\bar u \,i f''' \partial_x u \d x \,\le\,   
  \frac1{12} \int_\R |\partial_x u|^2 \d x + 3 \int_{|x|\le A}\beta^2 
  f'''^2 |u|^2 \d x  + 3 \int_{|x| \ge A} \beta^2 f'''^2 
  |u|^2 \d x ~.
\]
In the region $|x| \le A$, we have
\[
  3 \int_{|x| \le A}\beta^2 f'''(x)^2 |u|^2 \d x \,\le\, 
  3 \beta_0^2 K_3^2  \int_\R |u|^2 \d x \,\le\, 
  3 \beta_0^2 K_3^2  \int_\R (|\partial_x u|^2 + x^2 |u|^2)\d x~,
\]
and we choose $\beta_0 > 0$ small enough so that $3 \beta_0^2 K_3^2
\le 1/6$. In the region $|x| \ge A$ we know that $|f'''(x)| \le C_k 
|f'(x)|$, and in view of \eqref{betaprop} we find for $\epsilon$ 
small enough:
\[
  3 \int_{|x| \ge A}\beta^2 f'''(x)^2 |u|^2 \d x \,\le\,  
  \int_\R \frac{\beta}{12\eps} f'(x)^2 |u|^2 \d x~.
\]
Summarizing, we have shown that
\begin{align} \nonumber
 \frac{\D}{\D t}\,\Re\int_\R \beta (\partial_x \bar u)\,i f'(x)u \d x
 \,\le\, &-\int_\R \frac{3\beta}{4\epsilon}\,f'(x)^2 |u|^2 \d x
 + \frac13 \int_\R ( |\partial_x u|^2 + x^2 |u|^2) \d x \\ \label{id3'}
 &+ \int_\R \frac{5\alpha}{8} |\partial_x^2 u - x^2 u|^2 
  \d x + \int_\R \frac\gamma{2} f'(x)^2 |\partial_x u|^2 \d x~.
\end{align}

Finally we turn our attention to \eqref{id4}. This time the right-hand 
side contains only one additional term, which can be estimated 
as follows:
\begin{align*}
  -\Re \int_\R \gamma' f'(x)^2 \bar u \partial_x u \d x
  \,&\le\, \int_\R \frac \gamma4  f'(x)^2 |\partial_x u|^2 \d x +  
  \int_\R \frac{\gamma'^2}{\gamma} f'(x)^2 |u|^2 \d x \\
  \,&\le\, \int_\R \frac \gamma4  f'(x)^2 |\partial_x u|^2 \d x + 
  \int_\R \frac {\beta}{12\epsilon} f'(x)^2 |u|^2 \d x ~,
\end{align*}
where we have used the fact that $\gamma'^2/\gamma \ll \beta/\epsilon$
if $\epsilon$ is small, see \eqref{alphagamrel}, \eqref{betaprop}. 
It remains to bound the last term in the right-hand side of 
\eqref{id4}. Using a standard ($\epsilon$-independent) partition of 
unity and we decompose $\gamma$ as $\gamma = \gamma_1 + \gamma_2$, 
where~$\gamma_1$ and~$\gamma_2$ are nonnegative Lipschitz functions 
satisfying $\supp(\gamma_1) \subset \{|x| \le 2A\}$ and $\supp(\gamma_2) 
\subset \{|x| \ge 3A/2\}$. Then we observe that
\begin{align*}
  -2\Re \int_\R \gamma_1 f'(x) f''(x) \bar u\partial_x u\d x  
  \,&\le\, \int_\R \frac\beta{12\epsilon} f'(x)^2 |u|^2 \d x 
  + 12 K_2^2 \int_{|x| \le 2A} \frac{\epsilon\gamma^2}{\beta} 
  |\partial_x u|^2 \d x \\
  \,&\le\, \int_\R \frac\beta{12\epsilon} f'(x)^2 |u|^2 \d x 
  + \frac1{12}\int_\R |\partial_x u|^2 \d x~,
\end{align*}
where in the second inequality we have used the relation $\epsilon 
\gamma^2/\beta = 64 \beta^2$ and the fact that $\beta(x) \le \beta_0
2^{2k}$ if $|x| \le 2A$. We have also chosen $\beta_0 > 0$ small enough 
so that $768\,\beta_0^2\,K_2^2\,2^{4k} \le 1/12$. On the other hand, 
integrating by parts, we find
\[
  -2\Re \int_\R \gamma_2 f'(x) f''(x) \bar u\partial_x u\d x  \,=\, 
  \int_\R \gamma_2 (f''^2 + f' f''') |u|^2 \d x + 
  \int_\R \gamma'_2 f' f'' |u|^2 \d x~.
\]
We know that $|f''| + |f'''| \le C_k |f'|$ on the support of 
$\gamma_2$. Moreover, since $|\gamma'/\gamma| \le \frac32 
|\beta'/\beta| \le C$ by \eqref{betadef2}, it is clear that
$|\gamma_2'| \le C\gamma$ for some $C > 0$ independent of $\epsilon$. 
Thus
\[
  -2\Re \int_\R \gamma_2 f'(x) f''(x) \bar u\partial_x u\d x  
 \,\le\, C \int_\R \gamma f'(x)^2 |u|^2\d x \,\le\,  
 \int_\R \frac\beta{12\epsilon} f'(x)^2 |u|^2\d x~,
\]
because $\gamma \ll \beta/\epsilon$ if $\epsilon$ is small, see
\eqref{alphagamrel}, \eqref{betaprop}. We conclude that
\begin{align}\nonumber
  \frac12 \frac{\D}{\D t}\int_\R \gamma  f'^2|u|^2\d x 
  \,\le\, & -\int_\R \frac{3\gamma}4 f'(x)^2 (|\partial_x u|^2 + 
  x^2 |u|^2)\d x \\ \label{id4'} 
  & +\frac1{12}\int_\R |\partial_x u|^2 \d x + \int_\R \frac{\beta}
  {4\epsilon} f'(x)^2 |u|^2\d x~. 
\end{align}

Summarizing, if we collect estimates \eqref{id1}, \eqref{id2'},
\eqref{id3'}, and \eqref{id4'}, we obtain
\begin{align}\label{dotPhi}
  \Phi'(t) \,\le\, &-\frac14 \int_\R (|\partial_x u|^2 + x^2|u|^2)\d x 
  - \int_\R \frac{\beta}{4\epsilon} f'(x)^2 |u|^2 \d x \\ \nonumber
  &- \int_\R \frac{\alpha}8 |\partial_x^2 u - x^2 u|^2\d x 
  - \int_\R \frac{\gamma}4 f'(x)^2 (|\partial_x u|^2 + x^2 |u|^2)\d x~,
\end{align}
and as in Section~\ref{proofthm1} it is sufficient to keep only 
the first line in \eqref{dotPhi}. On the other hand, with our
choice of the function $\beta$, it is proved in Lemma~\ref{betamodif}
that there exists $M_0 > 0$ such that 
\begin{equation}\label{lowerbound1}
  \int_\R \Bigl(|\partial_x u|^2 + x^2 |u|^2 + \frac{\beta}{\epsilon}
  \,f'(x)^2 |u|^2\Bigr)\d x \,\ge\, \frac{M_0}{\eps^{\bar \nu}}
  \,\|u\|_{L^2}^2~, \quad \hbox{where} \quad \bar \nu \,=\, 
  \frac2{k+4} ~\cdotp
\end{equation}
Combining \eqref{lowerbound1} with the first line of \eqref{dotPhi}, 
we easily obtain 
\[
   \Phi'(t) \,\le\, -\frac18 \int_\R (|\partial_x u|^2 + x^2|u|^2)\d x 
   -\int_\R \frac{\beta}{8\epsilon} f'(x)^2 |u|^2 \d x
   -\frac{M_0}{8\epsilon^{\bar\nu}} \int_\R |u|^2 \d x~.
\]
Using now the upper bound \eqref{Philow}, we deduce that 
$\Phi'(t) \le -\eta \, \Phi(t)$ if
\[
  \eta \,=\, \min\left(\frac{1}{6\|\alpha\|_{L^\infty}}\,,\,
  \frac{1}{6\epsilon}\Bigl\|\frac\gamma\beta\Bigr\|_{L^\infty}^{-1}\,,\,
  \frac{M_0}{8\epsilon^{\bar \nu}}\right) \,=\,  \OO(\epsilon^{-\bar\nu})~.
\]
As in Section~\ref{proofthm1}, we conclude that $\Psi(\epsilon)^{-1}
\le C\,\epsilon^{\bar\nu} \log(2/\epsilon)$ for some $C > 0$ 
independent of $\epsilon$. Up to a logarithmic correction,
this proves the lower bound in \eqref{Psioptimal} with the rate 
$\epsilon^{\bar\nu}$, which is better than $\epsilon^{\nu}$ since 
$\nu=(k+2)^{-1} < \bar \nu$. 


\section{Resolvent estimates} \label{resolvent}

In this section we obtain precise estimates on the resolvent of
$H_\epsilon$ along the imaginary axis, and thereby prove
Theorem~\ref{thm2}. In doing so, we construct a nontrivial complex
domain $\omega_{\epsilon}$, parametrized by $\epsilon\in (0,1]$,
which avoids the pseudospectrum of $H_{\epsilon}$ as $\epsilon\to 0$ 
in the sense of Definition~\ref{de.pseudo}. The boundary of this
domain is a graph over the imaginary axis which exhibits a 
complicated structure involving various $\epsilon$-dependent scales, 
see Fig.~2 below. This phenomenon does not occur in familiar
examples such as self-adjoint operators, kinetic Fokker-Planck
operators, or complex harmonic oscillators, see \cite{DSZ}. 
Our proof will show that the value of the quantity~$\Psi(\epsilon)$ 
results from the competition between various microlocal models related 
either to the critical points of~$f$ or to the behavior of $f$ at infinity.

Throughout this section, we assume that $f : \R \to \R$ satisfies
Hypothesis~\ref{hypf}. In particular $\overline{f(\R)}= f(\R)\cup
\{0\}$ and $f$ has only a finite number of critical points. The finite 
set of critical values of $f$ is denoted by
\[
  \vc (f) \,=\, \Bigl\{f(x)\,;\, x \in \R\,,~f'(x) = 0\Bigr\}~.
\]
For any $\lambda \in \R$ and any $\epsilon \in (0,1)$, we define
\begin{equation}\label{eq.kappadef}
  \kappa(\epsilon,\lambda) \,=\, \|(H_{\epsilon}-i\lambda)^{-1}\|~. 
\end{equation}
The following proposition gives accurate bounds on $\kappa(\epsilon,
\lambda)$ in various parameter regimes: 

\begin{proposition}\label{pr.Hpseudoraff}
For $\epsilon \in (0,1)$ and $\lambda \in \R$, the quantity 
$\kappa(\epsilon,\lambda)$ defined in \eqref{eq.kappadef}
satisfies the following estimates: 
\begin{description}
\item[i)] If $\dist(\epsilon \lambda, f(\R))\ge \delta>0$, then
  $\kappa(\epsilon, \lambda)\le \epsilon/\delta$. 
\item[ii)] If $\dist(\epsilon\lambda, \vc (f)\cup \{0\})\ge 
 \delta > 0$, then $\kappa(\epsilon,\lambda)\le
 C_{\delta}\,\epsilon^{2/3}$. 
\item[iii)] If $\lambda=\lambda(\epsilon)$ is such that 
 $\lim_{\epsilon\to 0}\epsilon\lambda(\epsilon)=c\in \vc
 (f)\setminus \left\{0\right\}$, then
 $\limsup_{\epsilon\to 0}\epsilon^{-1/2}\kappa(\epsilon,\lambda
 (\epsilon)) \le C$. 
\item[iv)] For $\lambda=0$, the quantity $\kappa(\epsilon,0)$
 satisfies
\[
  \kappa(\epsilon,0) \,\le\, \left\{
  \begin{array}[c]{ll}
  C\,\epsilon^{\frac{2}{k+2}} & \text{if}\quad 0\not\in f(\R)~,\\
  C\,\epsilon^{\min\left\{\frac{2}{k+2},\frac{2}{3}\right\}}
  & \text{if}\quad 0\in f(\R)\setminus \vc (f)~,\\
  C \epsilon^{\min\left\{\frac{2}{k+2},\frac{1}{2}\right\}}
  &\text{if}\quad 0\in \vc (f)~.
\end{array}
\right.
\]
\item[v)] There exists $C > 1$ such that $\kappa(\epsilon,\lambda)
\le C\epsilon^{\frac{2}{k+4}}$ for all $(\epsilon,\lambda) \in (0,1)
\times \R$. Moreover, if $\kappa(\epsilon,\lambda)\ge C^{-1} 
\epsilon^{\frac{2}{k+4}}$, then $\lambda$ is comparable to 
$\epsilon^{-\frac{4}{k+4}}$. 
\end{description}
Finally all estimates in \textbf{i)}, \textbf{ii)}, \textbf{iii)}, 
\textbf{iv)}, and \textbf{v)} are optimal, in the sense that one can
find $\lambda = \lambda(\epsilon)$ so that the pair $(\epsilon,
\lambda(\epsilon))$ satisfies the required conditions as $\epsilon
\to 0$ and so that $\kappa(\epsilon,\lambda(\epsilon))$ is comparable 
to the upper bound in this limit.
\end{proposition}

Theorem~\ref{thm2} is of course a direct consequence of 
Proposition~\ref{pr.Hpseudoraff}, since
\[
  \Psi(\epsilon) \,=\, \min_{\lambda\in\R}\kappa(\epsilon,\lambda)^{-1}
  \,\in\, \left[M_4^{-1}\epsilon^{-\frac{2}{k+4}}\,,\,M_4
  \epsilon^{-\frac{2}{k+4}}\right]~, \quad \hbox{for some }
  M_4 \ge 1~.
\]
Proposition~\ref{pr.Hpseudoraff} also allows to localize the 
pseudospectrum of $H_\epsilon$ accurately: 

\begin{corollary}\label{co.pseudoloc}
The complex domain $\omega_\epsilon$ defined for $\epsilon\in (0,1]$ by 
\[
  \omega_{\epsilon} \,=\, \bigcup_{\lambda\in \R}\Bigl\{z\in \C\,;\,
  |z-i\lambda| \le \frac{1}{2\kappa(\epsilon,\lambda)}\Bigr\}
  ~\cup~ \Bigl\{z \in \C\,;\, \Re(z) \le 0\Bigr\}
\]
avoids the pseudospectrum of $H_{\epsilon}$ as $\epsilon\to 0$, 
whereas the domain
\[
  \tilde\omega_\epsilon \,=\, \{z \in \C\,;\,\Re(z) \le
  \Psi(\epsilon)\log(\epsilon^{-1})^2\} ~\cap~ (\C\setminus 
  \omega_\epsilon)
\]
meets the pseudospectrum of $H_{\epsilon}$ as $\epsilon\to 0$.
\end{corollary}

\proof If $|z-i\lambda| \le (2\kappa(\epsilon,\lambda))^{-1}$, the
resolvent formula shows that $\|(H_{\epsilon}-z)^{-1}\| \le
2\kappa(\epsilon, \lambda) \le C$. Moreover, we know that
$\|(H_{\epsilon}-z)^{-1}\| \le 1$ if $\Re(z) \le 0$, thus
$\omega_\epsilon$ avoids the pseudospectrum of $H_\epsilon$ as
$\epsilon\to 0$. Remark that $\omega_\epsilon$ contains the half-plane
$\{z \in \C\,;\, \Re(z) \le \Psi(\epsilon)/2\}$.  On the other hand,
if $\mu_\epsilon = \Psi(\epsilon)\log (\epsilon^{-1})^2$, it follows
from Proposition~\ref{pr.Hpseudoraff} that $\mu_\epsilon \gg
\Psi(\epsilon)(1+\log\Psi(\epsilon) + \log(\epsilon^{-1}))$
if~$\epsilon$ is small, hence the half-plane $\{\Re(z) \le
\mu_\epsilon\}$ meets the pseudospectrum of $H_\epsilon$ by
Lemma~\ref{pseudo}. The same conclusion holds of course for $\tilde
\omega_\epsilon$. \QED

\goodbreak
\begin{figure}[hbt]
\begin{center}
\vskip 0.8truecm{\baselineskip=0.8\baselineskip}
\includegraphics{fig2.epsi}
\end{center}
\vskip 0.8truecm{\baselineskip=0.8\baselineskip
\noindent \vbox{\noindent {\footnotesize {\bf Fig.~2:} The domain
  $\omega_\epsilon$ on the left-hand side of the solid curve avoids 
  the pseudospectrum of $H_{\epsilon}$ as $\epsilon\to 0$. The 
  picture on the right shows the geometry at the scale $\epsilon z 
  =\OO(1)$, while the left picture focuses on the region where 
  $\epsilon z$ is small. Here $R_{\epsilon} =  \{z \in \C\,;\, 
  \Re z\geq 0\,,\; \min f\le \epsilon \Im z\le \max f\}$ and  
  $cv(f)=\{c_1,c_2,c_3\}$.}}\par}
\end{figure}

\noindent To prove Proposition~\ref{pr.Hpseudoraff}, we start with a variant 
of the so-called IMS-localization formula, see~\cite{CFKS}. 
\begin{lemma}\label{le.IMS} 
Consider a Schr{\"o}dinger-type operator $Q = -\Delta + V$ in
$\R^d$, where $V$ is a measurable function. Take a locally 
finite partition of unity $\{\chi_j\}_{j \in J}$, where $\chi_{j}\in
C_0^{\infty}(\R^d,\R)$, such that
\begin{equation}\label{eq.L2cond}
  \sum_{j\in J}\chi_j(x)^2 \,=\, 1~, \quad \hbox{for all }
  x \in \R^d~,
\end{equation}
and 
\begin{equation}\label{eq.defm}
  m_1^2 \,\eqdef\, \sup_{x\in \R^d} \,\sum_{j\in J} |\nabla\chi_j(x)|^2
  \,<\, +\infty~, \qquad 
  m_2^2 \,\eqdef\, \sup_{x\in \R^d} \,\sum_{j\in J} (\Delta \chi_j(x))^2
  \,<\, +\infty~.
\end{equation}
Then the estimate
\begin{equation}
\label{eq.Qeqini}
2\|Qu\|^2 + 3m_2^2 \|u\|^2 + 8m_1^2 \|\nabla u\|^2
  \,\ge\, \sum_{j\in J} \|Q\chi_j u\|^2~,
\end{equation}
holds for any $u\in C_0^\infty(\R^d)$. Moreover, if $\Re V\ge 0$,
then
\begin{equation}
\label{eq.Qeq}
2\|Qu\|^2 + 3m_2^2 \|u\|^2 + 8m_1^2\Re\,\langle Qu\,,u\rangle
  \,\ge\, \sum_{j\in J} \|Q\chi_j u\|^2~.
\end{equation}
\end{lemma}

\proof
For any $\chi\in C_0^{\infty}(\R^d,\R)$ we have
\begin{eqnarray*}
  Q^*\chi^2 Q &=& Q^*\chi Q\chi+ Q^*\chi[\chi,Q]\\
  &=& \chi Q^*Q \chi + [Q^*,\chi]Q\chi + Q^*\chi[\chi,Q]\\
  &=& \chi Q^*Q \chi + [Q^*,\chi]\chi Q + Q^*\chi[\chi,Q]
    + [Q^*,\chi][Q,\chi]\\
  &=& \chi Q^*Q\chi - [\Delta,\chi]\chi Q + Q^*\chi[\Delta,\chi]
    + [\Delta,\chi][\Delta,\chi]~.
\end{eqnarray*}
Since
$$
  [\Delta,\chi] \,=\, 2(\nabla\chi)\cdot\nabla + (\Delta\chi)
  \,=\, 2\nabla\cdot(\nabla\chi) - (\Delta\chi)~,
$$
this implies
\[
  Q^*\chi^2 Q \,=\, \chi Q^* Q \chi - \nabla\cdot(\nabla\chi^2)Q
  + (\Delta \chi)\chi Q + Q^*(\nabla \chi^2)\cdot\nabla + 
  Q^*\chi(\Delta \chi) - R_\chi^*R_\chi~,
\]
with $R_\chi = 2(\nabla\chi)\cdot\nabla+(\Delta\chi)$. We now
apply this identity with $\chi=\chi_j$ and sum over $j \in J$. 
In view of \eqref{eq.L2cond}, the left-hand side reduces to 
$Q^*Q$, and the second and fourth terms in the right-hand side
disappear, so that 
\[
  Q^*Q \,=\, \sum_{j\in J}\chi_jQ^* Q\chi_j + \sum_{j\in J}
  \Bigl((\Delta\chi_j)\chi_j Q + Q^*\chi_j(\Delta\chi_j)\Bigr) 
  - \sum_{j\in J}R_{\chi_j}^* R_{\chi_j}~.
\]
Thus, for any $u\in C_0^\infty(\R^d,\R)$, we have
\begin{eqnarray*}
  \|Qu\|^2 &\geq& 
  \sum_{j\in J}\Bigl( \|Q\chi_j u\|^2 - (\|\chi_j Qu\|^2 + 
  \|(\Delta\chi_j)u\|^2) - (8\|(\nabla\chi_j)\cdot\nabla u\|^2
  + 2\|(\Delta \chi_j)u\|^2)\Bigr)\\
  &\geq & \Bigl(\sum_{j\in J} \|Q\chi_j u\|^2\Bigr)
  - \|Qu\|^2 - 3m_2^2\|u\|^2 -8m_1^2\|\nabla u\|^2~,
\end{eqnarray*}
which is \eqref{eq.Qeqini}. When $\Re V \geq 0$, the inequality
 $\|\nabla u\|^2 \le \Re\,\langle 
Q u\,,u\rangle$ implies \eqref{eq.Qeq}. \QED

\medskip
The idea is now to apply Lemma~\ref{le.IMS} to the one-dimensional
operator $H_\epsilon - i\lambda$, using a dyadic partition of unity.
This allows us to reduce a global problem on the whole real line to a
family of compactly supported problems, indexed by the parameter
$j \in \N$. The choice of a dyadic partition is convenient to take 
into account the precise behavior of the coefficients of $H_\epsilon$ 
as~$x \to \pm \infty$.

\begin{lemma}\label{le.reduc}
For $j\in \N$, $\epsilon > 0$, and $\lambda \in \R$, consider the 
operator
\begin{equation}\label{eq.Pdef}
  P_{j,\epsilon,\lambda} \,=\, -2^{-2j}\partial_x^2 + 2^{2j}x^2
  +\frac{i}{\epsilon}f(2^j x)-i\lambda~,
\end{equation}
and let
\begin{equation}\label{eq.Cjep}
  C_j(\epsilon,\lambda) \,=\, \inf\Bigl\{\|P_{j,\epsilon,\lambda}u\|
  \,;\, u\in C_0^\infty(\R)\,,~\supp u\subset K_j\,,~
  \|u\| = 1\Bigr\}~,
\end{equation}
where $K_0 = [-1,1]$ and $K_j = [-1,-1/4]\cup [1/4,1]$ for any 
$j > 0$. Then the quantity $\kappa(\epsilon,\lambda) = \|(H_\epsilon - 
i\lambda)^{-1}\|$ satisfies
\begin{equation}\label{eq.kappabdd}
  \Bigl(\inf_{j\in\N}C_j(\epsilon,\lambda)\Bigr)^{-1} \,\le\, 
  \kappa(\epsilon,\lambda) \,\le\, C \Bigl(\inf_{j\in\N}C_j
  (\epsilon,\lambda)\Bigr)^{-1}~,
\end{equation}
for some constant $C \ge 1$ independent of $\epsilon,\lambda$.
\end{lemma}

\begin{remark}\label{r.Cjeps}
It is clear that $C_j(\epsilon,\lambda) \ge 1$ for all $j \in \N$, 
$\epsilon > 0$, $\lambda \in \R$, because
\[
  \|u\|^2 \,\le\, \Re\,\langle P_{j,\epsilon,\lambda} u\,,u\rangle
  \,\le\, \|P_{j,\epsilon,\lambda} u\|\,\|u\|~,
\]
for all $u \in C_0^\infty(\R)$. 
\end{remark}

\proof
We choose a dyadic partition of unity $\{\chi_j\}_{j\in\N}$
such that
\[
  1 \,=\, \sum_{j=0}^{\infty}\chi_j(x)^2 \,=\, \chi_0(x)^{2}
  + \sum_{j=1}^{\infty}\tilde\chi\Bigl(\frac{x}{2^j}\Bigr)^2~,
\]
where $\chi_0, \tilde \chi \in C_0^\infty(\R)$ satisfy
\[
  \chi_0(x) \,=\, \left\{
  \begin{array}{rcl} 
  1 & \hbox{if} & |x| \le \frac34~, \\[1mm]
  0 & \hbox{if} & |x| \ge 1~,
  \end{array}\right. \qquad
  \tilde\chi(x) \,=\, \left\{
  \begin{array}{rcl} 
  1 & \hbox{if} & \frac12 \le |x| \le \frac34~, \\[1mm]
  0 & \hbox{if} & |x| \le \frac38 \hbox{ or } |x| \ge 1~.
  \end{array}\right.
\]
It is clear that such a partition exists and that the quantities 
$m_{1}^{2}$, $m_{2}^{2}$ defined by \eqref{eq.defm} are finite. 
Thus we can apply Lemma~\ref{le.IMS}  to the operator $Q = 
H_\epsilon - i\lambda = -\partial_x^2 + x^2 + \frac{i}{\epsilon}f(x)
-i\lambda$, for any~$\epsilon > 0$ and $\lambda \in \R$. Since 
$\|u\|^2 \le \Re\,\langle Qu\,, u\rangle \le \|Qu\|\,\|u\| \le 
\|Qu\|^2$ for all $u \in C_0^\infty(\R)$, it follows from 
\eqref{eq.Qeq} that
\[
  C^2 \|Qu\|^2 \,\ge\, \sum_{j=0}^\infty \|Q\chi_j u\|^2~,
\]
where $C^2 = 2 + 8m_1^2 + 3m_2^2$. Now, for any $j \in \N$, we define
\[
  v_j(x) \,=\, 2^{j/2}\chi_j(2^j x)u(2^j x)~, \quad 
  x \in \R~,
\]
so that $\supp v_j \subset \supp \chi_j(2^{j}\cdot)\subset K_j$
and $(P_{j,\epsilon,\lambda} v_j)(x) = 2^{j/2} (Q\chi_j u)(2^j x)$. 
If we denote $m(\epsilon,\lambda) = \inf_{j \in \N}C_j(\epsilon,\lambda)$, 
we thus find
\begin{equation}\label{eq.Qprelim}
  C^2 \|Q u\|^2 \,\ge\, \sum_{j=0}^\infty \|P_{j,\epsilon,\lambda}v_j\|^2
  \,\ge\, m(\epsilon,\lambda)^2 \sum_{j=0}^\infty \|v_j\|^2 
  \,\ge\, m(\epsilon,\lambda)^2 \|u\|^2~,
\end{equation}
because $\sum_j \|v_j\|^2 = \sum_j \|\chi_ j u\|^2 = \|u\|^2$. 
Since $Q = H_\epsilon -i\lambda$ and $C_0^\infty(\R)$ is dense in 
$L^2(\R)$, it follows from \eqref{eq.Qprelim} that $\kappa(\epsilon,
\lambda) = \|(H_\epsilon - i\lambda)^{-1}\| \le
C/m(\epsilon,\lambda)$, which is the upper bound in \eqref{eq.kappabdd}.

To prove the lower bound, we fix $\delta > 0$, $\epsilon > 0$, 
$\lambda\in \R$, and we take $j \in \N$, $v_j \in C_0^\infty(\R)$ 
such that $v_j \not\equiv 0$, $\supp v_j \subset K_j$, and 
$\|P_{j,\epsilon,\lambda}v_j\| \le (m(\epsilon,\lambda)+\delta)\|v_j\|$.
Setting $v_j(x) = 2^{j/2}u(2^j x)$, we find that $\|Qu\| \le 
(m(\epsilon,\lambda)+\delta)\|u\|$, hence $\kappa(\epsilon,\lambda) 
\ge (m(\epsilon,\lambda)+\delta)^{-1}$. Since $\delta > 0$ was 
arbitrary, this concludes the proof. \QED

\bigskip
\noindent\textbf{Proof of Proposition \ref{pr.Hpseudoraff}:}\null\\[2mm]
The proof of assertion \textbf{i)} is easy: if 
$\dist(\epsilon\lambda,f(\R)) \ge \delta$, then $\dist(i\lambda,
\Theta(H_\epsilon)) \ge \delta/\epsilon$, where $\Theta(H_\epsilon)$ 
is the numerical range defined in \eqref{Thetadef}. The last
inequality in \eqref{specbounds} then implies $\kappa(\epsilon,
\lambda) \le \epsilon/\delta$. 

In the remaining four cases, we start from Lemma~\ref{le.reduc} and
use \eqref{eq.kappabdd} to bound $\kappa(\epsilon,\lambda)$. 
It turns out to be convenient to rewrite the operator 
$P_{j,\epsilon,\lambda}$ in the equivalent form
\begin{equation}\label{eq.Pdef2}
  P_{j,\epsilon,\lambda} \,=\, \frac{1}{\epsilon 2^{kj}}
  \left[-\epsilon2^{(k-2)j}\partial_{x}^{2} + \epsilon
  2^{(k+2)j}x^{2}+i(2^{kj}f(2^{j}x)-2^{kj}\epsilon \lambda)\right]~,
\end{equation}
where $k > 0$ is the parameter that governs the asymptotic 
behavior of $f(x)$ as $|x| \to \infty$, see Hypothesis~\ref{hypf}.
For later use, we observe that
\begin{equation}\label{eq.Cfdef}
  C_f \,\eqdef\, \sup_{j\in\N} \,\sup_{x \in K_j}
  2^{kj}|f(2^j x)| \,<\, +\infty~. 
\end{equation}

\noindent\textbf{ii)} 
Suppose that $\dist(\epsilon\lambda,\vc(f)\cup\{0\})\ge \delta$. 
Without loss of generality, we also assume that $\epsilon|\lambda|
\le \|f\|_{L^\infty} + \delta$, because otherwise we can use 
the estimate established in \textbf{i)}. For any $u \in C_0^\infty(\R)$ 
with $\supp u\in K_j$ and $u \not\equiv 0$, we have the lower bound
\[
  \frac{\|P_{j,\epsilon,\lambda}u\|}{\|u\|} \,\ge\, \frac{|\Im\,
  \langle P_{j,\epsilon,\lambda}u\,,u\rangle|}{\|u\|^2} \,=\, 
  \frac{|\langle [2^{kj}f(2^j\cdot)-2^{kj}\epsilon\lambda]u
  \,,u\rangle|}{\epsilon 2^{kj}\|u\|^2} \,\ge\, 
  \frac{1}{\epsilon}\Bigl(\epsilon|\lambda| - \frac{C_f}{2^{kj}}
  \Bigr)~,
\]
where $C_f$ is given by \eqref{eq.Cfdef}. Since $\epsilon|\lambda| \ge
\delta$, we deduce that $C_j(\epsilon,\lambda)\ge \delta/(2\epsilon)$
whenever $j$ is large enough so that $2^{kj} \ge 2C_f/\delta$. Thus
only a finite number of indices $j$ have to be considered, and the
problem is therefore reduced to finding a lower bound on  
the quantity $\|(H_{\epsilon}-i\lambda)u\|$ when $u\in C_0^\infty
(\{x \in \R\,;\, |x| < R_{\delta}\})$, for some $R_\delta > 0$. 
On a bounded domain, we can drop the bounded term $x^2$ in 
$H_\epsilon$ and only consider the operator $Q = -\partial_x^2 + 
\frac{i}{\epsilon}(f(x)-\epsilon\lambda)$. Take two functions 
$\theta_0$, $\theta_1 \in C^\infty(\R,\R)$ such that $\supp\theta_0
\subset[-2,2]$, $\supp \theta_1\subset(-\infty,-1] \cup [1,+\infty)$, 
and $\theta_0(x)^2 + \theta_1(x)^2 = 1$ for all $x \in \R$. Given 
$\sigma>0$ which will be fixed below, we consider the new
partition of unity
\begin{equation}\label{eq.chi+-}
  \chi_0(x)^2 + \chi_1(x)^2 \,\equiv\, 1~, \quad\hbox{with}\quad 
  \chi_j(x) \,=\, \theta_j\Bigl(\frac{f(x)-\epsilon\lambda}
  {\epsilon^{\sigma}}\Bigr) \quad \hbox{for}\quad j = 0,1~.
\end{equation}
The quantities \eqref{eq.defm} for this partition satisfy 
$m_1^2\le C_1 \epsilon^{-2\sigma}$ and $m_2^2\le C_1 
\epsilon^{-4\sigma}$ for some $C_1 > 0$, hence \eqref{eq.Qeq}
implies
\begin{equation}\label{eq.Qtot}
  2\|Qu\|^2 + \frac{3C_1}{\epsilon^{4\sigma}}\|u\|^2
 + \frac{8C_1}{\epsilon^{2\sigma}}\|Qu\|\,\|u\|
 \,\ge\, \|Q(\chi_0 u)\|^2 + \|Q(\chi_1 u)\|^2~.
\end{equation}
The last term is easily bounded from below using the support 
condition on $\chi_1$, which yields
\begin{equation}\label{eq.Q1}
  \|\chi_1 u\| \,\|Q(\chi_1 u)\| \,\ge\, \frac{1}{\epsilon}
  |\langle (f(x)-\lambda \epsilon)\chi_1 u\,,\chi_1 u\rangle|
  \,\ge\, \epsilon^{\sigma-1} \|\chi_1 u\|^2~.
\end{equation}
To bound $Q(\chi_0 u)$, we observe that $f^{-1}(\epsilon\lambda)$
is a finite set $\{x_1,\dots,x_n\}$ because $f$ has a only a finite
number of critical points. Since $\dist(\epsilon\lambda,\vc(f)\cup
\{0\})\ge \delta$, it follows that for $\epsilon > 0$ small enough
the support of $\chi_0$ is the disjoint union of $n$ intervals 
$I_1,\dots,I_n$ of size $ \OO(\epsilon^\sigma)$ centered at the points 
$x_1,\dots,x_n$. In particular, $\chi_0 u = \sum_{\ell=1}^n u_\ell$ 
with $\supp u_\ell\cap \supp u_{\ell'} = \emptyset$ when $\ell \neq 
\ell'$. Inside $I_\ell$, the operator $Q$ is well approximated
by $Q_\ell = -\partial_x^2+\frac{i}{\epsilon}f'(x_\ell)(x-x_\ell)$, 
because
\begin{eqnarray}\nonumber
  \|Q u_\ell\|^2 &\ge& \frac{1}{2} \|Q_\ell u_\ell\|^2 -
  \left\|\frac{f(x)-f(x_\ell)-f'(x_\ell)(x-x_\ell)}{\epsilon}\,u_\ell
  \right\|^2 \\\label{eq.Quell}
  &\ge& \frac{1}{2}\|Q_\ell u_\ell\|^2 - \frac{C_2\,
  \|f''\|_{L^\infty}^2\,\epsilon^{4\sigma}}{\epsilon^{2}}
  \,\|u_\ell\|^2~.
\end{eqnarray}
On the other hand, the operator $Q_\ell$ is unitarily equivalent to
the microlocal model
\begin{equation}\label{eq.micromod}
  \tilde{Q}_\gamma \,=\, \gamma^{2/3}(-\partial_y^2 \pm iy)~, \quad 
  \hbox{where} \quad \gamma \,=\, \frac{|f'(x_\ell)|}{\epsilon}~,
\end{equation}
which satisfies $\|\tilde{Q}_{\gamma} u\| \ge C \gamma^{2/3}\|u\|$. 
Actually, the operator $P = -\partial_y^2 \pm iy$ has a compact
resolvent, since 
\[
  \|Pu\|^2  \,=\, \|\partial_y^2 u\|^2 + \|yu\|^2 \pm 2i\langle 
  u\,,\,\partial_y u\rangle \,\ge\, \frac12 \|\partial_y^2 u\|^2
  + \|yu\|^2 - C\|u\|^2~,
\]
with an empty kernel ($Pu = 0$ implies $\|\partial_y u\|^2 =
\Re\langle Pu\,,\, u\rangle = 0$), and it is therefore invertible.  We
refer to Chap.~27 in \cite{Hor} or to \cite{DSZ} for a geometric
analysis of similar and more general models. Alternatively, one can
use the Lie-algebra approach developed in \cite{RoSt}, \cite{HeNo} if
one sets $-\partial_y^2 +y\partial_t = -X_1^2+X_0$ with $X_1 =
\partial_y$, $X_0 = y\partial_t $ and if one considers the value
$\tau=1$ for the frequency variable $\tau$ dual to $t\in \R$.
Summarizing, we have shown that
\begin{equation}\label{eq.Q0}
  \|Q(\chi_0 u)\|^2 \,=\, \sum_{\ell=1}^n \|Q u_\ell\|^2 
  \,\ge\, \Bigl(C_3 \epsilon^{-4/3} - C_2 \|f''\|_{L^\infty}^2
  \,\epsilon^{4\sigma - 2}\Bigr) \|\chi_0 u\|^2~,
\end{equation}
for some $C_3 > 0$. We now assume that $\sigma > 1/6$, so that 
$\epsilon^{4\sigma-2} \ll \epsilon^{-4/3}$ if $\epsilon$ is small. 
Replacing \eqref{eq.Q1}, \eqref{eq.Q0} into \eqref{eq.Qtot}, we thus 
find
\[
  C_4 \|Qu\|^2 + \frac{C_4}{\epsilon^{4\sigma}}\|u\|^{2} \,\ge\, 
  C_4^{-1}\min\left\{\epsilon^{2\sigma-2},\epsilon^{-4/3}\right\}
  \,\|u\|^2~,
\]
for some $C_4 > 0$. Finally, we suppose that $\sigma < 1/3$, so that
$\epsilon^{-4\sigma} \ll \epsilon^{-4/3} \ll \epsilon^{2\sigma-2}$. We
thus obtain the estimate $2C_4^2 \|Qu\|^2 \ge \epsilon^{-4/3}\|u\|^2$, 
which proves that $\kappa(\epsilon,\lambda) \le C \epsilon^{2/3}$. 

\medskip
\noindent\textbf{iii)} The assumption $\lim_{\epsilon\to 0}\epsilon
\lambda = c\in\vc(f)\setminus\{0\}$ implies that $\epsilon|\lambda|
\ge \delta$ for some fixed $\delta > 0$ if~$\epsilon$ is small enough.
Thus we can reduce the analysis to a bounded spatial domain like in
case~\textbf{ii)}. By assumption~$f^{-1}(c)$ is a finite set which contains at
least one critical point of $f$, but in general this set contains
noncritical points as well. Using a partition of unity, we can treat
the noncritical points separately and estimate their contributions as
in case \textbf{ii)}. So, for simplicity, we assume here that $f^{-1}(c)$
consists of critical points only. We shall consider two different
regimes, depending on how fast $\epsilon \lambda$ converges to $c$ as
$\epsilon \to 0$, and then check after some iterations that they
provide a complete information for the general assumption
$\lim_{\epsilon\to 0}\epsilon\lambda=c$.

\noindent\textbf{a)} We first study the case when $\epsilon\lambda$ 
converges slowly to $c$. More precisely, we assume that  
\begin{equation}\label{eq.sigma12}
  \epsilon^{\sigma_1} \,\le\, |\epsilon\lambda-c| \,\le\, 
  \epsilon^{\sigma_2}~,
\end{equation}
where $0 \le \sigma_2 < \sigma_1 < 1/2$ and $3\sigma_2 > 5\sigma_1-1$.
If $\sigma_2 = 0$, we also suppose that $\epsilon\lambda \to c$ as
$\epsilon \to 0$. Several intervals $[\sigma_{2},\sigma_{1}]$ will be
fixed iteratively in step \textbf{c)} below. Our goal is to show 
that $\|Q u\| \ge C\epsilon^{-1/2}\|u\|$ for all $u \in
C_0^\infty(\R)$, where $Q = -\partial_x^2 + \frac{i}{\epsilon}
(f(x)-\epsilon\lambda)$. We choose $\sigma \in (\sigma_1,1/2)$ 
such that
\begin{equation}\label{eq.sigmacond}
  \frac{2\sigma_1}{3} + \frac16 \,<\, \sigma \,<\, \frac13 + 
  \frac{\sigma_1}{3} - \frac{\sigma_1-\sigma_2}{2}~,
\end{equation}
and we use again the partition of unity defined by \eqref{eq.chi+-}.
If $x\in \supp \chi_0$ or $x\in \supp\chi_1'$, we have $|f(x)-
\epsilon\lambda|\le 2\epsilon^\sigma$ and using \eqref{eq.sigma12}
we easily find
\[
  \frac{\epsilon^{\sigma_1/2}}{C_2} \,\le\, \frac{|f(x)-c|^{1/2}}
  {C_1} \,\le\, |f'(x)| \,\le\, C_1 |f(x)-c|^{1/2} \,\le\,
  C_2\,\epsilon^{\sigma_2/2}~,
\]
for some $C_1, C_2 > 0$. These estimates allow to give precise 
bounds on the quantities $m_1^2$ and $m_2^2$ defined in
\eqref{eq.defm}. Since
\begin{eqnarray*}
 \chi_j'(x) &=& \frac{f'(x)}{\epsilon^\sigma}\,\theta_j'
 \Bigl(\frac{f(x)-\lambda\epsilon}{\epsilon^\sigma}\Bigr)~,\\
 \chi_j''(x) &=& \frac{f''(x)}{\epsilon^\sigma}\theta_j'\,\Bigl(
 \frac{f(x)-\lambda\epsilon}{\epsilon^\sigma}\Bigr)
 + \Bigl(\frac{f'(x)}{\epsilon^\sigma}\Bigr)^2 \,\theta_j''
 \Bigl(\frac{f(x)-\epsilon\lambda}{\epsilon^\sigma}\Bigr)~,
 \quad j = 0,1~,
\end{eqnarray*}
we obtain
\[
  m_1^2 \,\le\, C_3 \,\epsilon^{\sigma_2-2\sigma}~, \quad 
  \hbox{and}\quad m_2^2 \,\le\, C_3^2 (\epsilon^{-2\sigma}+
  \epsilon^{2\sigma_2-4\sigma}) \,\le\, 2 C_3^2 \,\epsilon^{2\sigma_2-
  4\sigma}~,
\]
for some $C_3 > 0$, and it follows from \eqref{eq.Qeq} that
\begin{equation}\label{eq.llll}
 3\|Qu\|^2 + 22 C_3^2\, \epsilon^{2\sigma_2-4\sigma} \|u\|^2 \,\ge\, 
 \|Q(\chi_0 u)\|^2 + \|Q(\chi_1 u)\|^2~.
\end{equation}
By \eqref{eq.Q1} we have $\|Q(\chi_1 u)\| \ge \epsilon^{\sigma-1}
\|\chi_1 u\|$, and to bound $\|Q(\chi_0 u)\|$ we proceed as in 
case \textbf{ii)}. Denoting $f^{-1}(\epsilon\lambda) = \{x_1,\dots,x_n\}$, 
we decompose as before $\chi_0 u = \sum_{\ell=1}^n u_\ell$, and 
we observe that every $x \in \supp u_\ell$ satisfies $|x-x_\ell|
\le 2 C_2\,\epsilon^{\sigma-\sigma_1/2}$, because
\[
  2\epsilon^\sigma \,\ge\, |f(x)-\epsilon\lambda| \,\ge\, 
  |x-x_\ell| \,\inf \Bigl\{|f'(y)|\,;\, y \in \supp u_\ell\Bigr\}  
  \,\ge\, C_2^{-1}\,\epsilon^{\sigma_1/2}|x-x_\ell|~.
\]
As in \eqref{eq.Quell}, we thus find
\[
  \|Qu_\ell\|^2 \,\ge\, \Bigl(\frac{C |f'(x_\ell)|^{4/3}}{\epsilon^{4/3}}
  -\frac{4 C_2^4 \|f''\|_{L^\infty}^2 \epsilon^{4\sigma-2\sigma_1}}
  {\epsilon^2}\Bigr)\|u_\ell\|^2 \,\ge\, C_4 
  \,\epsilon^{\frac{2\sigma_1}{3} -\frac43}\|u_\ell\|^2~,
\]
for some $C_4 > 0$, because $\epsilon^{-4/3}|f'(x_\ell)|^{4/3} \ge
C_2^{-4/3} \epsilon^{\frac{2\sigma_1}{3} -\frac43} \gg
\epsilon^{4\sigma-2\sigma_1-2}$ by \eqref{eq.sigmacond}. Summing
over~$\ell$ as in \eqref{eq.Q0}, we obtain the desired lower bound on
$\|Q(\chi_0 u)\|^2$, and returning to \eqref{eq.llll} we arrive at
\[
  C_5 \|Qu\|^2 + C_5\,\epsilon^{2\sigma_2-4\sigma} \|u\|^2 \,\ge\, 
  C_5^{-1}\min\left\{\epsilon^{2\sigma-2}\,,\,\epsilon^{
  \frac{2\sigma_1}{3} -\frac43}\right\}\|u\|^2~,
\]
for some $C_5 > 0$. Since $\epsilon^{2\sigma_2-4\sigma} \ll 
\epsilon^{\frac{2\sigma_1}{3} -\frac43} \ll \epsilon^{2\sigma-2}$ 
by \eqref{eq.sigmacond}, we conclude that
\[
   2 C_5^2 \|Qu\|^2 \,\ge\, \epsilon^{\frac{2\sigma_1}{3} -\frac43}
   \|u\|^2 \,\ge\, \epsilon^{-1}\|u\|^2~, \quad \hbox{for all }
   u \in C_0^\infty(\R)~.
\]
\noindent\textbf{b)} We now assume that $|\epsilon \lambda - c| \le 
\epsilon^\sigma$, for some $\sigma \in (\frac13,\frac12)$, and 
we use again the partition of unity defined by \eqref{eq.chi+-}. 
As before, we have $|f'(x)|\le C_1|f(x)-c|^{1/2}\le C_2\,
\epsilon^{\sigma/2}$ for all $x\in \supp\chi_0$ and all $x 
\in \supp \chi_1'$, and it follows that the quantities \eqref{eq.defm} 
satisfy $m_1^2\le C_3\,\epsilon^{-\sigma}$, $m_2^2\le C_3^2\,
\epsilon^{-2\sigma}$. Thus \eqref{eq.Qeq} becomes
\[
  3\|Qu\|^2 + 22 C_3^2\,\epsilon^{-2\sigma} \|u\|^2 \,\ge\, 
  \|Q(\chi_0 u)\|^2 + \|Q(\chi_1 u)\|^2 \,\ge\, \|Q(\chi_0 u)\|^2
  + \epsilon^{2\sigma - 2}\|\chi_1 u\|^2~.
\]
Since $\epsilon^{-2\sigma} \ll \epsilon^{-1} \ll \epsilon^{2\sigma-2}$, 
it is sufficient to show that $\|Q(\chi_0 u)\|^2 \ge C\epsilon^{-1}
\|\chi_0 u\|^2$. To this end, we consider the set $f^{-1}(c) = 
\{x_1,\dots,x_n\}$ and we decompose as before $\chi_0 u = 
\sum_{\ell=1}^n u_\ell$. In the support of $u_\ell$, we can
use the quadratic approximation 
\[
  Q_\ell \,=\, -\partial_x^2 + \frac{i}{\epsilon}\Bigl[
  \frac{f''(x_\ell)}{2}(x-x_\ell)^{2} - (\epsilon\lambda-c)\Bigr]~.
\]
Indeed, if $x\in\supp u_\ell$, then $|x-x_\ell|\le  C_4 
\epsilon^{\sigma/2}$ for some $C_4 > 0$, and it follows that
\[
  \|Qu_\ell\|^2 \,\ge\, \frac{1}{2}\|Q_\ell u_\ell\|^2
  - \frac{C_{4}^{6}\|f'''\|_{L^\infty}^2 \epsilon^{3\sigma}}
  {36 \epsilon^2} \|u_\ell\|^2 \,=\,
  \frac{1}{2}\|Q_\ell u_\ell\|^2 - C_5 \,\epsilon^{3\sigma-2}
  \|u_\ell\|^2~,
\]
for some $C_5 > 0$. On the other hand, the operator $Q_\ell$ is 
unitarily equivalent to the microlocal model
\[
  \gamma^{1/2}(-\partial_y^2 \pm i y^2-i\mu)~, \quad 
  \hbox{where} \quad \gamma \,=\, \frac{|f''(x_\ell)|}{2\epsilon}
  \quad \hbox{and} \quad \mu \gamma^{1/2} \,=\, \frac{\epsilon\lambda -c}
  {\epsilon}~.
\]
Using the methods presented in \cite{Hor} (Chap.~27), \cite{HeNo}, 
\cite{DSZ}, or even reproducing the analysis of \textbf{ii)}  leading to
\eqref{eq.micromod}, we find  
\[
  \|(-\partial_y^2 \pm i y^2-i\mu)v\| \,\ge\, C(1+|\mu|^{1/3})
  \|v\| \,\ge\, C \|v\|~, \quad \hbox{for all } v\in C_0^\infty(\R)~.
\]
Since $\epsilon^{3\sigma-2} \ll \epsilon^{-1}$, this shows that
$\|Qu_\ell\|^2 \ge C\,\epsilon^{-1}\|u_\ell\|^2$, hence also 
$\|Q(\chi_0 u)\|^2 \ge C\,\epsilon^{-1}\|\chi_0 u\|^2$. 
As in step \textbf{a)}, we conclude that $\|Q u\|^2 \ge 
C\,\epsilon^{-1}\|u\|^2$ for all $u \in C_0^\infty(\R)$. 

\noindent\textbf{c)} Take $\sigma_1^0 = \frac{11}{30}\in
(\frac{1}{3},\frac{1}{2})$ and, for any $n\in \N$, let $\sigma_2^n =
\sigma_1^{n+1} = \frac{11}{6}\sigma_1^n - \frac{1}{3}$. Since 
$\frac{11}{6} > 1$ and since the solution to $x = \frac{11}{6}x 
-\frac{1}{3}$ is $x = \frac{2}{5} > \sigma_1^0$, it is clear that
$\sigma_2^n < \sigma_1^n$ for all $n$ and $\sigma_2^n \to -\infty$ as 
$n \to \infty$.  Let $n_0 \in \N$ be the smallest integer for which 
$\sigma_2^n \le 0$. As $\frac{11}{6} > \frac{5}{3}$, the condition 
$3\sigma_2^n > 5\sigma_1^n-1$ is satisfied for all $n \le n_0$. 
Applying step \textbf{a)} to all intervals $[\max\{0,\sigma_2^n\},
\sigma_1^n]$ for $n = 0,\dots,n_0$, we obtain the lower bound
$\|(H_\epsilon-i\lambda)u\| \ge C\,\epsilon^{-1/2}\|u\|$ whenever
$\lambda$ satisfies \eqref{eq.sigma12} with $[\sigma_2,\sigma_1] 
\subset [0,\frac{11}{30}]$. In other words, there exists a constant 
$K > 0$ such that
\[
  \left(
  \begin{array}[c]{c}
   |\epsilon\lambda(\epsilon)-c| \ge \epsilon^{\frac{11}{30}}\\ 
  \lim_{\epsilon\to 0}\epsilon\lambda(\epsilon)=c
  \end{array}\right) \quad\Rightarrow\quad
  \left(\limsup_{\epsilon\to 0}\epsilon^{-1/2}\kappa(\epsilon,
  \lambda(\epsilon))\le K\right)~.
\]
Meanwhile step \textbf{b)} provides the same conclusion (possibly 
with another constant $K$) under the assumption $|\epsilon\lambda
(\epsilon)-c|\le \epsilon^{\frac{11}{30}}$. This concludes the 
proof of case \textbf{iii)}. 

\medskip
\noindent\textbf{iv)} We consider the operator $P_{j,\epsilon,\lambda}$
given by \eqref{eq.Pdef2} for $\lambda = 0$ and $j \ge 1$. For 
any $u \in C_0^\infty(\R)$ with $\supp u \subset K_j = \{\frac14 \le 
|x| \le 1\}$, we have
\begin{align*}
 &\|u\|\,\|P_{j,\epsilon,0}u\| \,\ge\, |\Re\,\langle P_{j,\epsilon,0}u
 \,,u \rangle| \,\ge\, 2^{2j}\int_{K_j} x^2 |u(x)|^2\d x \,\ge\, 
 2^{2j-4} \|u\|^2~, \\[1mm]
 &\|u\|\,\|P_{j,\epsilon,0}u\| \,\ge\, |\Im\,\langle P_{j,\epsilon,0}u
 \,,u \rangle| \,\ge\, \frac{1}{\epsilon 2^{kj}}\int_{K_j}
 2^{kj}|f(2^j x)| |u|^2 \d x \,\ge\, \frac{m_j}{\epsilon 2^{kj}}
 \,\|u\|^2~,
\end{align*}
where $m_j = \inf\{2^{kj}|f(2^j x)|\,;\,\frac14 \le |x| \le 1\}$. 
From Hypothesis~\ref{hypf} we know that $m_j \to 1$ as $j \to \infty$, 
so we can take $J \in \N$ large enough so that $m_j \ge 1/2$ whenever
$j \ge J$. Thus
\[
  C_{j}(\epsilon,\lambda) \,\ge\, \frac12 \Bigl(2^{2j-4} 
  + \frac{1}{\epsilon 2^{kj+1}}\Bigr) \,\ge\, 
  C\,\epsilon^{-\frac{2}{k+2}}~, \quad \hbox{for all } j \ge J~. 
\]
The case when $j \in \{0,\dots,J\}$ corresponds to a bounded spatial
domain, and can be treated exactly as in \textbf{ii)} and
\textbf{iii)}.  We find that $\|H_\epsilon u\| \ge
C\epsilon^{-\sigma}\|u\|$ where $\sigma = 1$ if $f^{-1}(0) =
\emptyset$, $\sigma = 2/3$ if $0 \in f(\R)\setminus\vc(f)$, and
$\sigma = 1/2$ if $0 \in \vc(f)$.

\medskip
\noindent\textbf{v)} We want to show that $\kappa(\epsilon,\lambda)
\le C_1\,\epsilon^{\frac{2}{k+4}}$ for all $(\epsilon,\lambda) \in 
(0,1) \times \R$. In view of Lemma~\ref{le.reduc}, we have to verify 
that 
\begin{equation}\label{eq.new1}
  C_j(\epsilon,\lambda) \,\ge\, C_2\,\epsilon^{-\frac{2}{k+4}}~, 
  \quad \hbox{for all } j \in \N \,\hbox{ and all }\, 
  (\epsilon,\lambda) \in (0,1) \times \R~.
\end{equation}
As was already observed, the analysis of $C_j(\epsilon,\lambda)$ for
$j$ in a finite set $\{0,\dots,J\}$ can be performed as in
\textbf{ii)}, \textbf{iii)} and yields (in the worst case) the lower
bound $C_j(\epsilon,\lambda) \ge C_J\,\epsilon^{-1/2}$, which is
already better than \eqref{eq.new1}.  Thus is it sufficient to
consider large values of $j$. Fix such a $j$ and choose $u \in
C_0^\infty(\R)$ such that $\supp u\subset K_j = \{\frac14 \le |x| \le
1\}$, $\|u\| = 1$, and $\|P_{j,\epsilon, \lambda}u\| \le
2C_j(\epsilon,\lambda)$. Proceeding as in \textbf{iv)}, we easily find
\begin{equation}\label{eq.new2}
  \|P_{j,\epsilon,\lambda}u\| \,\ge\, 2^{2j-4}~, \quad
  \hbox{and} \quad \|P_{j,\epsilon,\lambda}u\| \,\ge\,  
  \frac{\inf_{x \in K_j}|g_j(x)|}{\epsilon 2^{kj}}~, 
\end{equation}
where $g_j(x) = 2^{kj}f(2^jx) - 2^{kj}\epsilon\lambda$. If $2^j \ge 
\epsilon^{-\frac{1}{k+4}}$, the first inequality already implies 
\eqref{eq.new1}. In the converse case, we integrate by parts and
obtain the relation
\begin{equation}\label{eq.new3}
  \|P_{j,\epsilon,\lambda}u\|^2 + 2 \|u\|^2 \,=\, 
  \|Q_{j,\epsilon,\lambda}u\|^2 + 2 \|x u_x\|^2 + 2^{4j}\|x^2 u\|^2
  \,\ge\, \|Q_{j,\epsilon,\lambda}u\|^2~,
\end{equation}
where $Q_{j,\epsilon,\lambda} = P_{j,\epsilon,\lambda} - 2^{2j}x^2$.
Thus $\|P_{j,\epsilon,\lambda}u\| \ge \|Q_{j,\epsilon,\lambda}u\| 
- \sqrt{2}\|u\|$, and combining this estimate with \eqref{eq.new2}
we obtain the sharper result
\begin{equation}\label{eq.new4}
  2C_j(\epsilon,\lambda) \,\ge\, \|P_{j,\epsilon,\lambda}u\| 
  \,\ge\, \frac13 \Bigl(2^{2j-4} + \frac{\inf_{x \in K_j}|g_j(x)|}
  {\epsilon 2^{kj}} + \|Q_{j,\epsilon,\lambda}u\| -\sqrt{2}\Bigr)~. 
\end{equation}
We next observe that
\[
  Q_{j,\epsilon,\lambda} \,=\, \frac{1}{\epsilon 2^{kj}}\,
  (-h^2\partial_x^2 + ig_j)~, \quad \hbox{where} \quad 
  h^2 \,=\, \epsilon 2^{(k-2)j} \,=\,  \OO(\epsilon^{\frac{6}{k+4}})~.
\]
From Hypothesis~\ref{hypf}, we know that the derivative $g_j'(x) = 
2^{(k+1)j}f'(2^j x)$ converges to $-kx|x|^{-k}$ as $j \to \infty$, 
uniformly on $K_j$. In particular, there exists $C_3 > 0$ such that
$|g_j'(x)| \ge C_3$ and $|g_j''(x)| \le C_3^{-1}$ for all $x \in K_j$ 
if $j$ is sufficiently large. Thus, using a partition of unity and 
a Taylor expansion like in case \textbf{ii)}, we can again reduce the 
analysis to the microlocal model \eqref{eq.micromod}. We obtain
the following lower bound
\[
  \|(-h^2\partial_x^2 + ig_j)u\| \,\ge\, C_4 \,h^{2/3}\|u\|~,
  \quad \hbox{for all } u \in C_0^\infty(\R) \hbox{ with }
  \supp u \subset K_j~.
\]
Returning to \eqref{eq.new4}, we arrive at the estimate
\begin{equation}\label{eq.new5}
  2 C_j(\epsilon,\lambda) \,\ge\, \frac13 \Bigl(2^{2j-4} +
  \frac{C_4\,\epsilon^{1/3} \,2^{(k-2)j/3}}{\epsilon 2^{kj}}
  - \sqrt{2}\Bigr) \,\ge\, C_5\,\epsilon^{-\frac{2}{k+4}}~,
\end{equation}
which proves \eqref{eq.new1}. 

Assume now that there is a family $(\epsilon,\lambda)$ with 
$\lambda = \lambda(\epsilon)$ such that $\epsilon \to 0$ and 
$\kappa(\epsilon,\lambda) \ge C_6 \,\epsilon^{\frac{2}{k+4}}$ 
for some $C_6 > 0$. By Lemma~\ref{le.reduc}, there exists $C_7 > 0$ 
such that, for all $(\epsilon,\lambda)$, 
\begin{equation}\label{eq.Cjcond}
  C_j(\epsilon,\lambda) \,\le\, C_7\,\epsilon^{-\frac{2}{k+4}}~, 
  \quad \hbox{for some } j = j(\epsilon,\lambda)\in \N~.
\end{equation}
It is clear that $j(\epsilon,\lambda) \to \infty$ as $\epsilon 
\to 0$, because if $j$ stays in a finite set $\{0,\dots,J\}$ we
have the lower bound $C_j(\epsilon,\lambda) \ge C_J\,\epsilon^{-1/2}$
which is incompatible with \eqref{eq.Cjcond}. We also know 
that $2^j \le C\,\epsilon^{-1/(k+4)}$, otherwise \eqref{eq.new4}
would contradict \eqref{eq.Cjcond} if $\epsilon$ is small. Thus 
we can use estimate \eqref{eq.new5} which implies, in view of   
\eqref{eq.Cjcond},
\begin{equation}\label{eq.new6}
  C_8^{-1}\,\epsilon^{-\frac{1}{k+4}} \,\le\, 2^j \,\le\, 
  C_8\,\epsilon^{-\frac{1}{k+4}}~, \quad \hbox{for some }
  C_8 \ge 1~.
\end{equation}
Next we choose $u$ as before and we deduce from \eqref{eq.new2}, 
\eqref{eq.Cjcond} that
\[
  2C_7\,\epsilon^{-\frac{2}{k+4}} \,\ge\, 2C_j(\epsilon,\lambda)
  \,\ge\, \|P_{j,\epsilon,\lambda}u\| 
  \,\ge\, \frac{\inf_{x \in K_j}|g_j(x)|}{\epsilon 2^{kj}}~\cdotp
\]
As $g_j(x) = 2^{kj}f(2^jx) - 2^{kj}\epsilon\lambda$, this implies
by direct calculation
\begin{equation}\label{eq.new7}
  \frac{1}{\epsilon 2^{kj}}\,\inf_{x\in K_j}2^{kj}f(2^j x) - 
  2C_7\,\epsilon^{-\frac{2}{k+4}} \,\le\, \lambda \,\le\, 
  \frac{1}{\epsilon 2^{kj}}\,\sup_{x\in K_j}2^{kj}f(2^j x) +
  2C_7\,\epsilon^{-\frac{2}{k+4}}~.
\end{equation}
Since $2^{kj}f(2^j x) \to |x|^{-k}$ uniformly in $K_j$ as $j \to \infty$ 
by Hypothesis~\ref{hypf}, and since $2^j$ satisfies \eqref{eq.new6}, 
it follows immediately from \eqref{eq.new7} that
\[
  C_9^{-1}\,\epsilon^{-\frac{4}{k+4}} \,\le\, \lambda \,\le\,  
  C_9\,\epsilon^{-\frac{4}{k+4}}~, 
\]
for some $C_9 \ge 1$, which is the desired estimate. 

Finally, it remains to check the optimality of the various 
estimates established so far. This is a rather straightforward 
task, because the arguments we gave to bound $\kappa(\epsilon,
\lambda)$ from above also indicate how to choose an appropriate
test function $u \in C_0^\infty(\R)$ which ``saturates'' the
inequality and gives the corresponding lower bound. We shall 
provide the details only in the last case \textbf{v)}, because 
this implies the upper bound on $\Psi(\epsilon)$ in 
Theorem~\ref{thm2}. 

Fix $x_0 \in (1/4,1)$ and take $j \in \N$ very large. We define 
$\epsilon> 0$, $\lambda > 0$, and $h > 0$ by the following relations
\begin{equation}\label{eq.epshdef}
  2^j \,=\, \epsilon^{-\frac{1}{k+4}}~, \qquad h^2 \,=\, \epsilon 
  2^{(k-2)j}~, \qquad \epsilon\lambda \,=\, f(2^j x_0)~.
\end{equation}
Next, we choose $v \in C_0^\infty(\R)$ such that $\|v\| = 1$ and
$\supp v \subset [-1,1]$, and we denote
\begin{equation}\label{eq.uhdef}
  u_h(x) \,=\, \frac{1}{h^{1/3}}\,v\Bigl(\frac{x-x_0}{h^{2/3}}
  \Bigr)~, \quad x \in \R~.
\end{equation}
It is clear that $u_h \in C_0^\infty(\R)$, $\|u_h\| = 1$, and that 
$\supp u_h \subset K_j$ if $h > 0$ is sufficiently small, that is, 
if $j$ is sufficiently large. Recalling that
\[  
  P_{j,\epsilon,\lambda} \,=\, \frac{1}{\epsilon 2^{kj}}
  \Bigl(-h^2 \partial_x^2 + h^{2/3} x^2+ ig_j(x)\Bigl)~, \qquad
  \hbox{where} \quad g_j(x) \,=\, 2^{kj}f(2^j x)- 2^{kj}
  \epsilon \lambda~,
\]
we claim that there exists $C_0 > 0$ independent of $j$ (hence
of $\epsilon,\lambda$) such that
\begin{equation}\label{eq.Pupper}
  \|P_{j,\epsilon,\lambda}u_h\| \,\le\, C_0\,\frac{h^{2/3}}{\epsilon
  2^{kj}} \,=\, C_0 \,\epsilon^{-\frac{2}{k+4}}~.
\end{equation}
This implies that $C_j(\epsilon,\lambda) \le C_0 
\,\epsilon^{-\frac{2}{k+4}}$, hence $\kappa(\epsilon,\lambda)
\ge C_0^{-1}\,\epsilon^{\frac{2}{k+4}}$ by \eqref{eq.kappabdd}. 

Checking \eqref{eq.Pupper} is of course straightforward. 
First, using \eqref{eq.uhdef}, we find $\|h^2 \partial_x^2 u_h\|
= h^{2/3}\|v''\|$. Next, since $x^2 \le x_0^2 + 2|x-x_0|$ for all
$x \in K_j$, we have $\|x^2 u_h\| \le x_0^2 \|v\| + 2h^{2/3}
\|v\| \le C$. Finally, since $g_j(x_0) = 0$ by our 
choice of $\lambda$, we have for all $x \in [1/4,1]$
\[
  |g_j(x)| \,\le\, |x-x_0|\,\sup_{z \in [1/4,1]} |g_j'(z)| \,\le\, 
  C |x-x_0|~,
\]
where the constant $C$ does not depend on $j$ by Hypothesis~\ref{hypf}. 
Thus $\|g_j u_h\| \le Ch^{2/3} \|v\|$, and the proof of 
\eqref{eq.Pupper} is complete. \QED

\figurewithtex 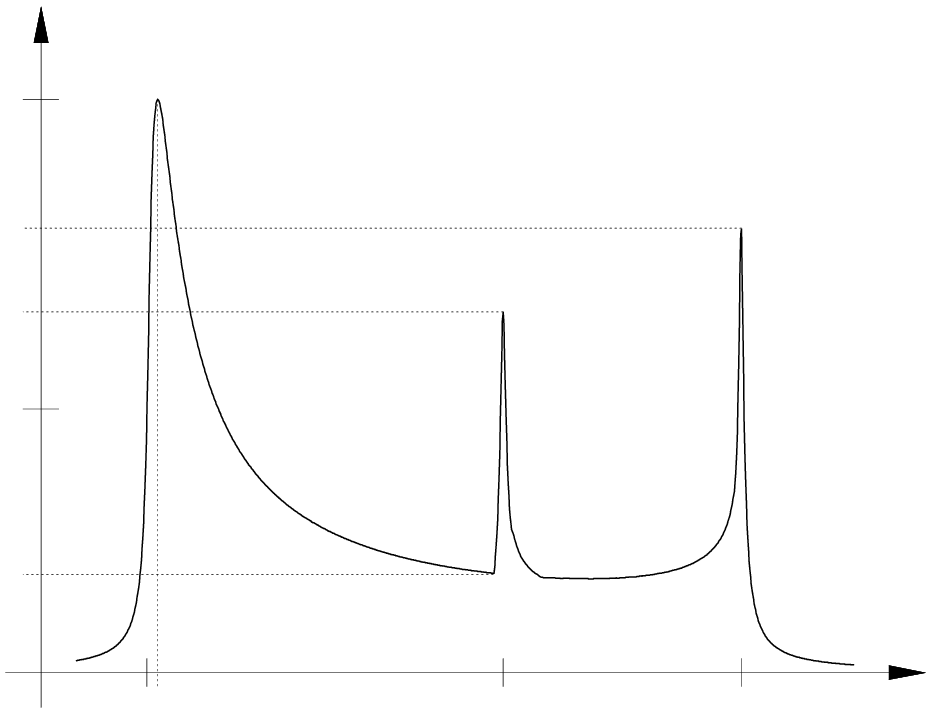 Fig3.tex 7.500 10.500 {\bf Fig.~3:} The
resolvent $\kappa(\epsilon,\lambda) =
\|(H_{\epsilon}-i\lambda)^{-1}\|$ is represented for $\epsilon =
2^{-14}$ and $f$ given by \eqref{eq.fex}.\cr

To illustrate the various regimes appearing in 
Proposition~\ref{pr.Hpseudoraff}, we computed numerically the
function $\kappa(\epsilon,\lambda)$ in the particular case
when
\begin{equation}\label{eq.fex}
  f(x) \,=\, \frac{1+3x^2}{(1+x^2/3)^3}~, \quad x \in \R~.
\end{equation}
Up to a multiplicative factor, this function satifies
Hypothesis~\ref{hypf}, with $k = 4$.  It has exactly two critical
values $c_1 = 1$, $c_2 = 27/16$, which correspond to the critical
points $x_1 = 0$ and $x_2 = \pm 1$ respectively. The result for
$\epsilon = 2^{-14}$ is displayed in Fig.~3. The large peak near the
origin corresponds to the regime described in \textbf{v)} above, while
the two sharp spikes are associated to the critical values of $f$, as
explained in \textbf{iii)}. The predictions of \textbf{i)},
\textbf{ii)} are also confirmed, but the value $\epsilon = 2^{-14}$ is
not small enough in this example to test the prediction of
\textbf{iv)}: we know from Proposition~\ref{pr.Hpseudoraff} that
$\kappa(\epsilon,0) = \OO(\epsilon^{1/3})$ since $k = 4$ and $0 \notin
f(\R)$, but in Fig.~3 the value computed for $\lambda = 0$ is still
clearly below the $ \OO(\epsilon^{1/2})$ spikes associated to the
critical values of $f$.  


\section{Spectral lower bounds}\label{spectral}

This section is devoted to the detailed study of the particular
example \eqref{fex}. We shall first prove Proposition~\ref{pr.speclb},
and then illustrate the optimality of the lower bound
\eqref{Sigmaconj} using numerical calculations and heuristic
semiclassical arguments.  As was mentioned in the introduction, the
estimate \eqref{Sigmaconj} implies that $\|e^{-tH_{\epsilon}}\| \le
C\,e^{-\mu_\epsilon t}$ for some $C > 0$ if $\mu_\epsilon = \frac12
M_5\, \epsilon^{-\nu'}$, but it follows from Lemma~\ref{elem} that $C
\ge e^{-1}\exp(\mu_\epsilon/ \Psi(\epsilon))$ where
$\mu_\epsilon/\Psi(\epsilon) = \OO(\epsilon^{\bar\nu-\nu'}) \to
\infty$ as $\epsilon \to 0$. It is very likely that such an
exponential bound in the semiclassical framework cannot be obtained
while remaining in the $C^{\infty}$ category. Thus, like in the study
of resonances (see \cite{AgCo,CFKS,HeSj,HiSi}), we shall exploit the
analyticity properties of the function \eqref{fex} by performing a
complex deformation. In this example we shall observe an interesting
phenomenon when comparing the pseudospectral and spectral estimates.
In contrast to what happens for resonances, the large peak of the
quantity $\kappa(\epsilon,\lambda) = \|(H_{\epsilon}-
i\lambda)^{-1}\|$ in the region $\lambda \sim \epsilon^{-4/(k+4)}$
(see Fig.~3) is not due to the presence of an eigenvalue nearby, for
the deformed operator.

\bigskip\noindent\textbf{Proof of Proposition~\ref{pr.speclb}:}
We shall prove the lower bound \eqref{Sigmaconj} in two steps. 
We first observe that, in view of Proposition~\ref{pr.Hpseudoraff}, 
it is sufficient to show that there exist positive constants 
$c_1$, $c_2$ such that
\begin{equation}\label{sector}
  \sigma(H_{\epsilon}) \,\cap\, \Bigl\{z\in \C\,;\, c_1 \Re z \le 
  |\Im z| \le \frac{c_2}{\epsilon}\Bigr\} \,=\, \emptyset~,
\end{equation}
for $\epsilon > 0$ small enough. Indeed, assume for the moment that
\eqref{sector} holds. If $\lambda \in \R$ is such that $|\lambda| 
\ge c_2\,\epsilon^{-1}$, the proof of Proposition~\ref{pr.Hpseudoraff} 
shows that $\kappa(\epsilon,\lambda) = \|(H_\epsilon-i\lambda)^{-1}\| 
\le C \epsilon^{1/2}$ for some $C > 0$, hence any $z \in 
\sigma(H_\epsilon)$ with $|\Im z| \ge c_2\,\epsilon^{-1}$ necessarily 
satisfies $\Re z \ge C^{-1}\epsilon^{-1/2}$. On the other hand, the 
spectrum of $H_\epsilon$ in the strip $|\Im z| \le c_2\,\epsilon^{-1}$
is concentrated in the sector $|\Im z| \le c_1 \Re z$ by
\eqref{sector}, hence to prove the lower bound \eqref{Sigmaconj} it 
is sufficient to consider the smaller strip $|\Im z| \le 
\epsilon^{-\nu'}$. We thus take $\lambda \in \R$ with $|\lambda| \le
\epsilon^{-\nu'}$ and we estimate $\kappa(\epsilon,\lambda)$ 
as in the proof of Proposition~\ref{pr.Hpseudoraff}. Using exactly 
the same notations, we find for all $j \in \N$:
\begin{equation}\label{auxcj}
  C_j(\epsilon,\lambda) \,\ge\, \frac12\Bigl(2^{2j-4} + 
  \frac{1}{\epsilon 2^{kj}}\,\inf_{x \in K_j} |2^{kj}f(2^jx) - 
  2^{kj}\epsilon\lambda|\Bigr)~,
\end{equation}
see e.g. \eqref{eq.new2}. As usual it is sufficient to consider large 
values of $j$, in which case $2^{kj}f(2^jx) \approx |x|^{-k}$ for all
$x \in K_j = \{x \in \R\,;\, \frac14 \le |x| \le 1\}$. If $2^{kj}
\epsilon|\lambda| \le 1/4$, the infimum in \eqref{auxcj} is larger 
than $1/2$ for large $j$, hence
\[
  C_j(\epsilon,\lambda) \,\ge\, \frac12\Bigl(2^{2j-4} + 
  \frac{1}{\epsilon 2^{kj+1}}\Bigr) \,\ge\, C\,\epsilon^{-
  \frac{2}{k+2}} \,\ge\, C\,\epsilon^{-\nu'}~.
\]
If $2^{kj} \epsilon|\lambda| \ge 1/4$ and $|\lambda| \le
\epsilon^{-\nu'}$, then $2^{kj} \ge C \epsilon^{\nu'-1}$ hence
$C_j(\epsilon,\lambda) \ge 2^{2j-5} \ge C \epsilon^{-\nu'}$. Thus we
have shown that $\kappa(\epsilon,\lambda)^{-1} \ge C \inf_j C_j
(\epsilon,\lambda) \ge C \epsilon^{-\nu'}$ if $|\lambda| \le
\epsilon^{-\nu'}$, from which we deduce that any $z \in
\sigma(H_\epsilon)$ with $|\Im z| \le \epsilon^{-\nu'}$ necessarily
satisfies $\Re z \ge C^{-1}\epsilon^{-\nu'}$. Summarizing, we have
shown that there exists $M_5 > 0$ such that $H_\epsilon$ has no
spectrum in the half-plane $\{\Re z \le M_5 \epsilon^{-\nu'}\}$
if~$\epsilon$ is sufficiently small, which is the desired result.

It remains to prove \eqref{sector}. Clearly this cannot be done using
only the pseudospectral estimates of Proposition~\ref{pr.Hpseudoraff},
because by Corollary~\ref{co.pseudoloc} the domain $\{z \in \C\,;\,
c_1 \Re z \le |\Im z| \le c_2\,\epsilon^{-1}\}$ meets the
pseudospectrum of $H_\epsilon$ as $\epsilon \to 0$. Instead we shall
implement a complex deformation method using the group of dilations
$(U_\theta\varphi)(x) = e^{\theta/2}\varphi(e^\theta x)$, which are
unitary operators when $\theta\in \R$. If~$f$ is given by \eqref{fex},
the multiplication operator by $(i/\epsilon)f(x)$ is a dilation
analytic perturbation of the harmonic oscillator hamiltonian $H_\infty
= -\partial_x^2 + x^2$, according to Definition~8.1 in \cite{CFKS}.
This implies that
\[
  H_\epsilon(\theta) \,=\, U_\theta H_\epsilon U_\theta^{-1} \,=\, 
  -e^{-2\theta}\partial_x^2 + e^{2\theta}x^2 + \frac{i}{\epsilon 
  (1+e^{2\theta}x^2)^{k/2}}
\]
defines an analytic family of type~(A) in the strip $S = \{\theta\in 
\C\,;\, |\Im \theta| < \pi/4\}$ (see~\cite{Ka}), with common domain $\DD =
D(H_{\infty})$. In particular the spectrum of $H_\epsilon(\theta)$,  
which is always discrete, does not depend on $\theta\in S$.

We choose $\theta= it_k$, where $t_k = \frac{\pi}{4(k+2)}$, and
we first observe that the operator $H_\epsilon(it_k)$ is still maximal 
accretive. Indeed, for all $(x,\xi) \in \R^2$, 
\begin{align*}
  e^{-2it_k}\xi^2 + e^{2it_k}x^2 \,&\in\, \Bigl\{z\in \C\,;\, 
  -2t_k \le \arg(z) \le 2t_k\Bigr\} \,\subset\, \{z \in \C\,;\,
  \Re(z) \ge 0\}~, \\
  \frac{i}{\epsilon}\,\frac{(1 + e^{-2it_k}x^2)^{k/2}}{|1+e^{2it_k}x^2|^{k}}
  \,&\in\, \Bigl\{z\in \C\,;\, \frac{\pi}{2}-kt_k \le \arg(z) \le
  \frac{\pi}{2}\Bigr\} \,\subset\, \{z \in \C\,;\,
  \Re(z) \ge 0\}~,
\end{align*}
hence $\Re \langle H(it_k)\varphi\,,\,\varphi\rangle \ge 0$ for all 
$\varphi\in \DD$. It is therefore sufficient to reproduce part of 
the analysis done in Proposition~\ref{pr.Hpseudoraff} to estimate
$\|(H_\epsilon-i\lambda)^{-1}\|$ for $\lambda \in \R$. Consider a 
partition of unity $\chi_0^2 + \chi_\infty^2 = 1$ with $\supp\chi_0
\subset (-1,1)$ and $\chi_0\equiv 1$ on $[-1/2,1/2]$. Applying 
the IMS localization formula \eqref{eq.Qeqini} to the operator 
$Q = e^{2it_k}(H_\epsilon(it_k)-i\lambda)$, and using the fact that
\[
  \|\partial_x u\|^2 \,\le\, \frac{1}{\cos(2t_k)} \Re
  \langle (H_\epsilon(it_k)-i\lambda) u\,,\,u\rangle
  \,\le\, \sqrt{2}\,\|u\|\,\|(H_\epsilon(it_k)-i\lambda)u\|~,
\] 
we obtain
\begin{equation}\label{eq.lowbdsp}
  3\|(H_\epsilon(it_k)-i\lambda)u\|^2 + C\|u\|^2 \,\ge\, 
  \|(H_\epsilon(it_k)-i\lambda)\chi_0 u\|^2 +
  \|(H_\epsilon(it_k)-i\lambda)\chi_\infty u\|^2~.
\end{equation}
We need to estimate both terms in the right-hand side. 

\noindent\textbf{i)} The localization of $u_0 =\chi_0 u$ in the 
interval $(-1,1)$ implies
\begin{align*}
  \|u_{0}\|\,\|(H_\epsilon(it_k)&-i\lambda)u_0\| 
  \,\ge\, \Im \langle(H_\epsilon(it_k)-i\lambda)u_0\,,\,u_0\rangle\\
  \,&\ge\, \frac{1}{\epsilon}\Re\,\langle (f(e^{it_k}x)-\lambda
  \epsilon)u_0\,,\, u_0\rangle -\sin(2t_k)\|\partial_x u_0\|^2\\
  \,&\ge\, \frac{1}{\epsilon}\inf_{|x| \le 1}\Bigl(\Re f(e^{it_k}x)
  -\epsilon\lambda\Bigr) \|u_0\|^{2} -\tan(2t_k)\|u_0\|\,
  \|(H_\epsilon(it_k)-i\lambda)u_0\|~.
\end{align*}
Since $|1+e^{2it_k}x^2|\le 2$ when $|x| \le 1$ and $-kt_k \le 
\arg(f(e^{it_k}x)) \le 0$, we have 
\[
  \Re f(e^{it_k}x) \,\ge\, \frac{\cos(kt_k)}{2^{k/2}} \,\ge\,  
  \frac{1}{2^{\frac{k+1}{2}}}~, \quad \hbox{for } |x| \le 1~.
\]
If we assume that $\epsilon|\lambda| \le c_2 \eqdef 2^{-\frac{k+2}{2}}$, 
we thus arrive at the lower bound
\[
  2\|(H_\epsilon(it_k)-i\lambda)\chi_0 u\| \,\ge\, 
  (1+\tan(2t_k))\|(H_\epsilon(it_k)-i\lambda)\chi_0 u\| \,\ge\, 
  c\,\epsilon^{-1}\|\chi_0 u\|~,
\]
where $c = 2^{-\frac{k+1}{2}}-2^{-\frac{k+2}{2}} > 0$.

\noindent\textbf{ii)} On the other hand, setting $u_\infty = 
\chi_\infty u$ and assuming that $u_\infty \neq 0$, we observe 
that
\[
  \frac{\|(H_\epsilon(it_k)-i\lambda)u_{\infty}\|}{\|u_\infty\|}
  \,\ge\, \left|\frac{\langle H_\epsilon(it_k) u_\infty\,,
  \,u_\infty\rangle}{\|u_\infty\|^2} - i\lambda\right| 
  \,\ge\, \inf_{z\in S_k}|z-i\lambda|~,
\]
where $S_k$ is any sector in the complex plane that is large 
enough to contain the quantity
\[
  \langle H_\epsilon(it_k) u_\infty\,,\,u_\infty\rangle \,=\,
  e^{-2it_k}\|\partial_x u_\infty\|^2 + e^{2it_k} \|xu_\infty\|^2 
  + \frac{e^{i(\pi/2-kt_k)}}{\epsilon}\Big\langle\frac{(e^{2it_k}+
  x^2)^{k/2}}{|1+e^{2it_k}x^2|^k}\,u_\infty\,,\, u_\infty\Big\rangle
\]
for any $u_\infty \in \DD$ with $\supp u_\infty \subset \R\setminus
[-\frac12,\frac12]$. If we denote
\[
  \delta_k \,=\, \frac{1}{2}\arg(e^{2it_k}+\frac{1}{4}) \in (0,t_k)~,
  \quad \hbox{so that}\quad k\delta_k \,\ge\, \max_{|x| \ge 1/2}
  \arg((e^{2it_k}+x^2)^{k/2})~,
\]
we have $\frac{\pi}{2}> \frac{\pi}{2}-k(t_k-\delta_k) > 
\frac{\pi}{2}-k t_k > 2t_k$, hence we can choose 
\[
  S_k  \,=\, \Bigl\{z\in \C\,;\, -2t_k \le \arg z \le \frac{\pi}{2}
  -k(t_k-\delta_k)\Bigr\}~,
\]
and there exists $c_k > 0$ such that $\inf_{z\in S_k}|z-i\lambda|
\ge c_k |\lambda|$.

If we now combine the lower bounds of both terms in the right-hand 
side of \eqref{eq.lowbdsp}, and if we use in addition the fact that 
$\|(H_\epsilon(it_k)-i\lambda)u\| \ge \cos(2t_k)\|u\|$, we conclude
that there exists $c_1 > 0$ such that
\begin{equation}\label{sectorial}
  \|(H_\epsilon(it_k)-i\lambda)u\| \,\ge\, C\,\min\Bigl\{
  \frac{c}{2\epsilon}\,,\,c_k |\lambda|\Bigr\} \|u\| \,\ge\, 
  \frac{2|\lambda|}{c_1}\,\|u\|~, \quad \hbox{for all } u \in \DD~,
\end{equation}
when $\lambda \in \R$ satisfies $|\lambda|\le c_2\,\epsilon^{-1}$.
Thus the deformed operator $H_\epsilon(it_k)$ is {\em sectorial}
in a neighborhood of the origin in the complex plane, uniformly 
for all $\epsilon \in (0,1]$. In particular, if $z = 
\mu + i\lambda$ with $0 < c_1 \mu \le |\lambda| \le c_2\,
\epsilon^{-1}$, it follows from \eqref{sectorial} that
\[
  \|(H_\epsilon(it_k)-z)^{-1}\| \,\le\, \frac{\|(H_\epsilon(it_k)-
  i\lambda)^{-1}\|}{1 - \mu\|(H_\epsilon(it_k)-i\lambda)^{-1}\|} 
  \,\le\, \frac{{\textstyle \frac12} c_1\,|\lambda|^{-1}}{1 - \mu(
  {\textstyle \frac12}c_1\,|\lambda|^{-1})} \,\le\, \frac{c_1}{|\lambda|}~,
\]
hence $z \notin \sigma(H_\epsilon(it_k)) \equiv \sigma(H_\epsilon)$. 
This proves \eqref{sector}, since we already know that $H_\epsilon$ 
has no spectrum in the half-plane $\Re z \le 0$. \QED

\bigskip Although Propositions~\ref{pr.Hpseudoraff} and
\ref{pr.speclb} provide the necessary information for comparing the
spectral bound $\Sigma(\epsilon)$ and the pseudospectral quantity
$\Psi(\epsilon)$, the analysis can be pushed further in the particular
example \eqref{fex}. As it is the case for $\Psi(\epsilon)$, the
actual behavior of $\Sigma(\epsilon)$ as $\epsilon \to 0$ results from
a competition between contributions due to the (unique) critical point
of $f$, and to the behavior of $f(x)$ as $|x| \to \infty$. Instead of
giving a complete proof, we briefly sketch the main idea, which is
adapted from the usual semiclassical techniques of harmonic
approximation (see for example \cite{Hel}) and complex deformation
(see for example \cite{Mar,Sj}). We also verify our conclusions using
numerical simulations.

We first write $\epsilon H_\epsilon  = -h^2 \partial_x^2 + V(x,h^2)$ 
where $h = \epsilon^{1/2}$ and $V(x,h^2) = h^2 x^2 + if(x)$, and 
we remark that $\epsilon H_\epsilon$ is not exactly a semiclassical
operator, because the potential $V$ depends on $h$ in a nontrivial
way which makes the comparison with quadratic approximations 
at critical points more difficult. In the sector $\{z\in \C\,;\,
|\arg z|\le \pi/4\}$, the function $V(z,\epsilon)$ has exactly 
three critical points: $z=0$, and $z=\pm z_{\epsilon}$, where
\[
  z_\epsilon \,=\, \Bigl\{\Bigl(\frac{ik}{2\epsilon}\Bigr)^{2\nu}
  -1\Bigr\}^{1/2} \,=\, \Bigl(\frac{ik}{2\epsilon}\Bigr)^\nu +
  \OO(\epsilon^\nu)~, \quad \hbox{where}\quad \nu = \frac{1}{k+2}~\cdotp
\]
The quadratic approximations of $H_\epsilon$ are respectively:

\medskip
\noindent\textbf{i)} At $z = 0$: $H_\epsilon^0 = -\partial_x^2 + 
\frac{i f(0)}{\epsilon} + (1+\frac{if''(0)}{2\epsilon})x^2$. 
The eigenvalues of $H_\epsilon^0$ are
\[
  \mu_n^0(\epsilon) \,=\, \frac{i}{\epsilon} + (2n+1)\omega_\epsilon~,
  \quad n \in \N~, \quad \hbox{where}\quad \omega_\epsilon \,=\,
  \Bigl(1 + \frac{i}{2\epsilon}f''(0)\Bigr)^{1/2}~.
\]

\noindent\textbf{ii)} At $z = \pm z_\epsilon$: $H_\epsilon^{\pm} = 
-\partial_x^2 + D_\epsilon + \Omega_\epsilon^2 x^2$, with
\begin{align*}
  D_\epsilon \,&=\, z_\epsilon^2 + \frac{i}{\epsilon}f(z_\epsilon)
  \,=\, \frac{k+2}{k}\,z_\epsilon^2 + \frac{2}{k} \,=\, \frac{k+2}{k}
  \Bigl(\frac{ik}{2\epsilon}\Bigr)^{2\nu} -1~, \\
  \Omega_\epsilon^2 \,&=\,  1 + \frac{i}{2\epsilon}f''(z_\epsilon)
  \,=\, (k+2)\,\frac{z_\epsilon^2}{1+z_\epsilon^2} \,=\, 
  (k+2)\Bigl(1 - \Bigl(\frac{2\epsilon}{ik}\Bigr)^{2\nu}\Bigr)~.
\end{align*}
\quad~ The eigenvalues of $H_\epsilon^{\pm}$ are $\nu_n^0(\epsilon)
= D_\epsilon + (2n+1)\Omega_\epsilon$, $n \in \N$.

\medskip\noindent After a complex rotation, these quadratic
approximations are transformed into real harmonic oscillators
$\mathbf{H}_\epsilon^0$ and $\mathbf{H}_\epsilon^\pm$, respectively.
The actual complex deformation has to be chosen carefully so that the
non-self-adjoint corrections to $\mathbf{H}_\epsilon^0$ (resp.
$\mathbf{H}_\epsilon^\pm$) still allow for a good control of the
resolvent norm of the deformed hamiltonian $\mathbf{H}_\epsilon$ when
$z$ is close to $\mu_n^0(\epsilon)$ (resp. $\nu_n^0(\epsilon)$). This
requires an accurate analysis, which is left for a future work.
In any case, on the basis of the arguments above, we expect that 
the full operator $H_\epsilon$ has two sets of eigenvalues 
$\{\mu_n(\epsilon)\}$ and $\{\nu_n(\epsilon)\}$ which satisfy:
\begin{align}\label{mundef}
  \mu_n(\epsilon) \,&=\, \frac{i}{\epsilon} + (2n+1)\omega_\epsilon
  + \OO(1)~, \quad n = 0,1,2,\dots\\ \label{nundef}
  \nu_n(\epsilon) \,&=\, D_\epsilon + (2n+1)\Omega_\epsilon + 
  \OO(\epsilon^\nu)~, \quad n = 0,1,2,\dots
\end{align}
This leads to the following conjecture: 

\begin{conjecture}\label{formal}
Fix $k > 0$ and let $f : \R \to \R$ be as in \eqref{fex}. 
Then the spectral bound of $H_\epsilon$ satisfies, when 
$\epsilon > 0$ is small, 
\begin{equation}\label{conjSigma}
  \Sigma(\epsilon) \,=\, \min\Bigl\{\Re(\mu_0(\epsilon))\,,\,
  \Re(\nu_0(\epsilon))\Bigr\}~,
\end{equation}
where $\mu_0(\epsilon)$ is given by \eqref{mundef} and
$\nu_0(\epsilon)$ by \eqref{nundef}.
\end{conjecture}

\noindent Conjecture~\ref{formal} predicts the following asymptotic 
expansions as $\epsilon \to 0$:
\[
  \Sigma(\epsilon) \,=\, \left\{ \begin{array}{lll}
  {\displaystyle \frac12 \Bigl(\frac{k}{\epsilon}\Bigr)^{1/2}}
  + \OO(1)~, & \hbox{if} & k \le 2~,\\[4mm]
  {\displaystyle \frac{k+2}{k}\Bigl(\frac{k}{2\epsilon}
  \Bigr)^{2\nu} \cos(\pi\nu) + \sqrt{k+2} - 1 + \OO(\epsilon^\nu)}
  & \hbox{if} & k > 2~, \end{array}\right.
\]
with $\nu = (k+2)^{-1}$, in full agreement with 
Proposition~\ref{pr.speclb}. 

\medskip To test these predictions, we have computed an approximation
of $\Sigma(\epsilon)$ using a standard finite difference scheme, for
three different choices of $k$ and for values of $\epsilon$ ranging
from $2^{-1}$ to $2^{-18}$. Calculations were made on a large interval
$[-L,L]$ with Dirichlet boundary conditions, using at least $500$ grid
points. The results are collected in Fig.~4, where the dots represent
the values computed numerically and the solid lines the theoretical
predictions, namely $\Re(\mu_0^0(\epsilon))$ for $k \le 2$ and
$\Re(\nu_0^0(\epsilon))$ for $k > 2$. The agreement is excellent and
leaves little doubt on the validity of Conjecture~\ref{formal}.

\figurewithtex 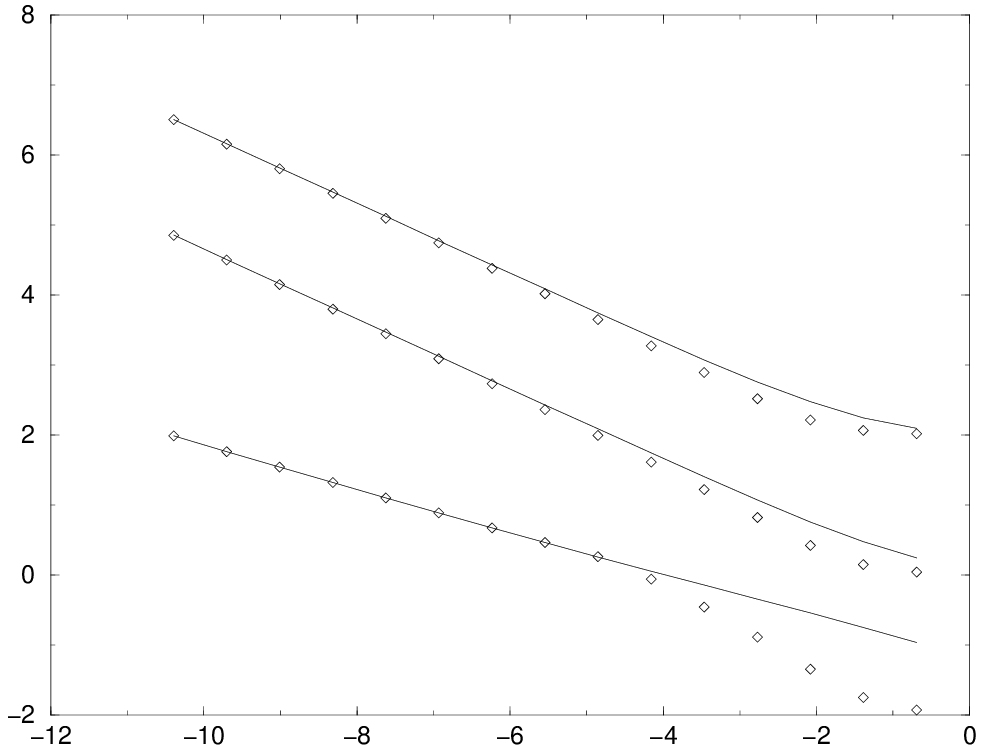 Fig4.tex 8.000 10.500
{\bf Fig.~4:} The logarithm of $\Sigma(\epsilon)$ is represented
as a function of $\log(\epsilon)$ for three choices of the parameter
$k$. The dots correspond to the values computed numerically, and the 
solid lines to the theoretical predictions given by 
$\Re(\mu_0^0(\epsilon))$ ($k = 1, 2$) and $\Re(\nu_0^0(\epsilon))$ 
($k = 4$). To improve readability, the upper curve has been 
shifted up by two units and the lower curve shifted down by two units. 
\cr

The heuristic arguments above give a correct prediction not
only for $\Sigma(\epsilon)$, but also for all eigenvalues of
$H_\epsilon$ with sufficiently low real part. Let us first consider
the case when $k = 4$ and $\epsilon = 2^{-18}$. In Table~1 below, we
have listed the five eigenvalues of $H_\epsilon$ with smallest real
part, together with the approximate eigenvalues $\mu_n^0(\epsilon)$
and $\nu_n^0(\epsilon)$ for $n = 0,\dots,4$. It is obvious that
$\nu_n^0(\epsilon)$ is an excellent approximation of
$\lambda_n(\epsilon)$.

\begin{center}
\begin{tabular*}{0.75\textwidth}{@{\extracolsep{\fill}}l c c c}
\hline
~$n$\rule[-3mm]{0mm}{8mm} 
 & $\lambda_n(\epsilon)$ & $\mu_n^0(\epsilon)$ & $\nu_n^0(\epsilon)$\\[1ex]
\hline
~0\rule[0mm]{0mm}{5mm} 
 & $106.18 + 60.48\,i$ & $512 + 2.616 \cdot 10^5\,i$ & $106.18 + 60.48\,i$ \\
~1 & $111.02 + 60.51\,i$ & $1536 + 2.606 \cdot 10^5\,i$ & $111.06 + 60.50\,i$ \\
~2 & $115.83 + 60.54\,i$ & $2560 + 2.596 \cdot 10^5\,i$ & $115.93 + 60.51\,i$ \\
~3 & $120.61 + 60.58\,i$ & $3584 + 2.586 \cdot 10^5\,i$ & $120.80 + 60.53\,i$ \\
~4 & $125.36 + 60.63\,i$ & $4608 + 2.575 \cdot 10^5\,i$ & $125.67 + 60.54\,i$ \\[1ex]
\hline
\end{tabular*}
\end{center}
\noindent{\footnotesize {\bf Table 1:}
The first few eigenvalues of $H_\epsilon$ for $k = 4$ and 
$\epsilon = 2^{-18}$, together with the predictions
$\mu_n^0(\epsilon)$ and $\nu_n^0(\epsilon)$.}

\bigskip
We next consider the case when $k = 2$, which is especially interesting 
because both $\Re(\mu_n^0(\epsilon))$ and $\Re(\nu_n^0(\epsilon))$ grow
with the same rate $\OO(\epsilon^{-1/2})$ as $\epsilon \to 0$. However
we see from \eqref{mundef}, \eqref{nundef} that $\Re(\mu_n^0(\epsilon)) 
\approx (2n+1)/(2\epsilon)^{1/2}$ while $\Re(\nu_n^0(\epsilon)) 
\approx (2/\epsilon)^{1/2} + 4n+1$. It follows that the eigenvalue
of $H_\epsilon$ with minimal real part is well approximated by 
$\mu_0^0(\epsilon)$, and thus corresponds to a semiclassical mode
near the origin, while the subsequent eigenvalues 
are well approximated by $\nu_n^0(\epsilon)$ for $n = 0,1,2,\dots$
and therefore correspond to semiclassical modes near infinity.
This prediction is clearly confirmed by the values listed in 
Table~2 below. It is also possible to verify that the eigenfunctions
corresponding to $\lambda_n(\epsilon)$ are nicely approximated 
either by an eigenfunction of the quadratic approximation 
$H_\epsilon^{0}$ at $z=0$ (if $n = 0$), or by an eigenfunction of the 
quadratic approximation $H_\epsilon^{\pm}$ at $z=\pm z_{\epsilon}$
(if $n = 1,2,\dots$). 

\begin{center}
\begin{tabular*}{0.75\textwidth}{@{\extracolsep{\fill}}l c c c}
\hline
~$n$\rule[-3mm]{0mm}{8mm}
& $\lambda_n(\epsilon)$ & $\mu_n^0(\epsilon)$ & $\nu_n^0(\epsilon)$\\[1ex]
\hline
~0\rule[0mm]{0mm}{5mm}  
& $44.54 + 4051\,i $ & $45.26 + 4050\,i  $ & $91.50 + 90.52\,i$ \\
~1 & $91.50 + 90.52\,i$ & $135.78 + 3960\,i $ & $95.48 + 90.54\,i$ \\
~2 & $95.47 + 90.54\,i$ & $226.30 + 3869\,i $ & $99.45 + 90.57\,i$ \\
~3 & $99.48 + 90.56\,i$ & $316.82 + 3779\,i $ & $103.43 + 90.59\,i$ \\
~4 & $103.2 + 90.58\,i$ & $407.34 + 3688\,i $ & $107.41 + 90.61\,i$ \\[1ex]
\hline
\end{tabular*}
\end{center}
\noindent{\footnotesize {\bf Table 2:}
The first few eigenvalues of $H_\epsilon$ for $k = 2$ and 
$\epsilon = 2^{-12}$, together with the predictions
$\mu_n^0(\epsilon)$ and $\nu_n^0(\epsilon)$.}

\bigskip Finally we turn our attention to the case when $k = 1$. Here
we expect that the eigenvalues of $H_\epsilon$ with lowest real part
are well approximated by $\mu_n^0(\epsilon)$, for $n = 0,1,2,\dots$.
This is indeed the case, although the spectrum of $H_\epsilon$ seems
more difficult to compute in this regime, due to numerical
instabilities. Nevertheless, for $\epsilon = 2^{-12}$, we were able to
check that at least the first three eigenvalues of $H_\epsilon$ are
correctly approximated by $\mu_n^0(\epsilon)$, see Table~3 below. The
errors are larger than in the preceding cases, but we should take into
account the fact that the remainder term in \eqref{mundef} is
$\OO(1)$.

\begin{center}
\begin{tabular*}{0.75\textwidth}{@{\extracolsep{\fill}}l c c c}
\hline
~$n$\rule[-3mm]{0mm}{8mm} 
& $\lambda_n(\epsilon)$ & $\mu_n^0(\epsilon)$ & $\nu_n^0(\epsilon)$\\[1ex]
\hline
~0\rule[0mm]{0mm}{5mm} 
& $31.46 + 4064\,i$ & $32.00 + 4064\,i $ & $242.63 + 419.00\,i$ \\
~1 & $93.34 + 4002\,i$ & $96.02 + 4000\,i $ & $246.09 + 419.01\,i$ \\
~2 & $153.2 + 3940\,i$ & $160.0 + 3936\,i $ & $249.55 + 419.01\,i$ \\[1ex]
\hline
\end{tabular*}
\end{center}
\noindent{\footnotesize {\bf Table 3:}
The first few eigenvalues of $H_\epsilon$ for $k = 1$ and 
$\epsilon = 2^{-12}$, together with the predictions
$\mu_n^0(\epsilon)$ and $\nu_n^0(\epsilon)$.}

\appendix

\section{Appendix}\label{appendix}

This appendix gathers the proofs of Lemmas~\ref{elem}, \ref{pseudo}, 
\ref{hatH}, and provides a variant of Lemma~\ref{hatH} which is needed in
Section~\ref{optimal}.

We first give a proof Lemma~\ref{elem}, which relates the quantities 
$\Sigma$ and $\Psi$ to the norm of the semigroup $e^{-t A}$ for a 
sectorial operator, as an application of the Laplace transformation.

\medskip\noindent
{\bf Proof of Lemma~\ref{elem}.} Let $A$ be a maximal accretive
operator in a Hilbert space $X$, with numerical range $\Theta(A)$ 
contained in the sector $\Delta_\alpha = \{z\in \C\,;\, |\arg z| \le 
\frac{\pi}{2}-2\alpha\}$ for some $\alpha \in (0,\frac{\pi}{4}]$. 

\noindent\textbf{i)} Assume first that $\|e^{-tA}\| \le C e^{-\mu t}$ 
for all $t \ge 0$. Then the resolvent
\[
  (A-z)^{-1} \,=\, \int_0^\infty e^{-tA}\,e^{tz}\d t
\] 
is defined (at least) for all $z \in \C$ with $\Re z < \mu$, hence
$\Sigma\ge \mu$. Moreover, taking $z = i\lambda \in i\R$ and using
the fact that $e^{-tA}$ is a semigroup of contractions in $X$, we
find
\[
  \|(A - i\lambda)^{-1}\| \,\le\, \int_0^\infty \|e^{-tA}\|\d t
  \,\le\, \int_0^\infty \min\{1\,,\,C e^{-\mu t}\}\d t \,=\, 
  \frac{1+\log(C)}{\mu}~\cdotp
\]
Taking the supremum over $\lambda \in \R$, we conclude that
$\Psi \ge \mu/(1+\log(C))$.

\noindent\textbf{ii)}
Conversely, to estimate the semigroup $e^{-tA}$ in terms of $\Sigma$ 
or $\Psi$, we use the inverse Laplace formula
\[
  e^{-tA} \,=\, \frac{1}{2\pi i} \int_{\Gamma(\mu,\alpha)}
  (A - z)^{-1} \,e^{-z t}\d z~,
\]
where $0 < \mu < \Sigma$ and $\Gamma(\mu,\alpha) = \Gamma_-(\mu,\alpha)
\cup \Gamma_0(\mu,\alpha) \cup \Gamma_+(\mu,\alpha)$ is the polygonal 
contour defined by 
\begin{align*}
  \Gamma_0(\mu,\alpha) \,&=\, \left\{z\in \C\,;\,\Re z = \mu\,,\, 
  |\arg z|\le \frac{\pi}{2}-\alpha\right\}~, \\
  \Gamma_\pm(\mu,\alpha) \,&=\, \left\{z\in \C\,;\,\Re z\ge \mu\,,\, 
  \arg z = \pm\Bigl(\frac{\pi}{2}-\alpha\Bigr)\right\}~,
\end{align*}
and oriented from $\Im z=-\infty$ to $\Im z= +\infty$. Note that 
$\Gamma(\mu,\alpha)$ lies entirely in the resolvent set of $A$ 
by construction. Since $\Re z = \mu$ when $z \in
\Gamma_0(\mu,\alpha)$, we easily obtain
\[
  \Bigl\| \int_{\Gamma_0(\mu,\alpha)} (A - z)^{-1} \,e^{-z t}\d z
  \Bigr\| \,\le\, N(A,\mu)\,\frac{2\mu}{\tan \alpha}\,e^{-\mu t}~,
\]
where $N(A,\mu) = \sup\{\|(A-z)^{-1}\|\,;\, \Re(z) = \mu\}$.
On the other hand, any $z \in \Gamma_+(\mu,\alpha)$ can be
parametrized as $z = x + (ix/\tan \alpha)$ with $x \ge \mu$, and the 
resolvent at this point can be estimated as follows: 
\[
  \|(A-z)^{-1}\| \,\le\, \frac{1}{\dist(z,\Theta(A))} \,\le\,  
  \frac{1}{\dist(z,\Delta_\alpha)} \,=\, \frac{1}{x}~.
\]
We thus find
\[
  \Bigl\| \int_{\Gamma_+(\mu,\alpha)} (A - z)^{-1} \,e^{-z t}\d z
  \Bigr\| \,\le\, \int_{\mu}^\infty \frac{1}{x}\, 
  e^{-tx}\,\frac{\d x}{\tan\alpha} \,\le\, \frac{1}{\tan\alpha}
  \,\frac{e^{-\mu t}}{\mu t}~,
\]
and the contribution of $\Gamma_-(\mu,\alpha)$ is estimated in 
exactly the same way. Collecting these bounds are using the 
fact that $\|e^{-tA}\| \le 1 \le 1/\tan \alpha$, we arrive at
\[
  \|e^{-tA}\| \,\le\, \frac{1}{\pi\tan\alpha}\left(
  \mu N(A,\mu)\,e^{-\mu t} + \min\Bigl\{\pi\,,\, \frac{e^{-\mu t}}
  {\mu t}\Bigr\}\right) \,\le\, \frac{1}{\pi\tan\alpha}
  \Bigl(\mu N(A,\mu) + 2\pi\Bigr)\,e^{-\mu t}~.
\]
\noindent\textbf{iii)} Assume that $0 < \mu < \Psi$. Using the
second resolvent formula and the definition of $\Psi$, we find for all 
$\lambda \in \R$:
\[
  \|(A-\mu-i\lambda)^{-1}\| \,\le\, \frac{\|(A-i\lambda)^{-1}\|}
  {1 - \mu\|(A-i\lambda)^{-1}\|} \,\le\, \frac{1}{\Psi - \mu}~,
\]
hence $N(A,\mu) \le (\Psi-\mu)^{-1}$. \QED

\medskip
We next give a proof of Lemma~\ref{pseudo}, which illustrates 
the pseudospectral nature of the quantity $\Psi(\epsilon)$. 

\medskip\noindent
{\bf Proof of Lemma~\ref{pseudo}.}\\[1mm]
\textbf{i)} If $\Re(z) \le 0$, we know that $\|(H_\epsilon - z)^{-1}\|
\le 1/\dist(z,\Theta(H_\epsilon)) \le 1$. If $0 < \Re(z) \le \kappa
\Psi(\epsilon)$ for some $\kappa \in (0,1)$, we have by 
Lemma~\ref{elem}-\textbf{iii)}:
\[
  \|(H_\epsilon-z)^{-1}\| \,\le\, N(H_\epsilon,\Re(z)) 
  \,\le\, \frac{1}{\Psi(\epsilon) - \Re(z)} \,\le\, 
  \frac{1}{\Psi(\epsilon)}\,\frac{1}{1-\kappa} \,\le\, 
  \frac{1}{1-\kappa}~\cdotp
\]
Thus $\|(H_\epsilon - z)^{-1}\|$ is uniformly bounded by
$(1-\kappa)^{-1}$ when $\Re(z) \le \kappa\Psi(\epsilon)$. 

\noindent\textbf{ii)} We argue by contraposition. Fix $K_0 \ge 1$, 
$N \in \N$, and assume that $\|(H_\epsilon - z)^{-1}\| \le K_0 
\,\epsilon^{-N}$ whenever $\Re(z) \le \mu_\epsilon$. We shall show
that there exists $K \ge 1$ (independent of $\epsilon$)
such that
\begin{equation}\label{contrapos}
  \mu_\epsilon \,\le\, K\Psi(\epsilon)\Bigl(1 + \log\Psi(\epsilon) 
  + \log(\epsilon^{-1})\Bigr)~, \quad \hbox{for all } \epsilon \in
  (0,1]~.
\end{equation}
Note that \eqref{contrapos} is automatically satisfied if
$\mu_\epsilon \le 1$, hence we assume from now on that $\mu_\epsilon 
\ge 1$. Applying Lemma~\ref{elem}-\textbf{ii)} with $\mu = \mu_\epsilon$,
and using the fact that $\alpha = \OO(\epsilon)$ when $A =
H_\epsilon$, we see that $\|e^{-tH_\epsilon}\| \le C(H_\epsilon,\mu_\epsilon)\, 
e^{-\mu_\epsilon t}$ for all $t \ge 0$, where
\[
  C(H_\epsilon,\mu_\epsilon) \,\le\, \frac{1}{\pi\tan\alpha}
  \Bigl(\mu_\epsilon K_0\,\epsilon^{-N} + 2\pi\Bigr) \,\le\, 
  \mu_\epsilon K_1 \,\epsilon^{-N-1}~, 
\]
for some $K_1 \ge 1$. Next, using the lower bound on $\Psi$
given by Lemma~\ref{elem}-\textbf{i)}, we obtain
\[
  \mu_\epsilon \,\le\, \Psi(\epsilon)\Bigl(1 + \log C(H_\epsilon,
  \mu_\epsilon)\Bigr) \,\le\, \Psi(\epsilon)\Bigl(1 + \log(\mu_\epsilon)
  + \log(K_1\,\epsilon^{-N-1})\Bigr)~.
\]
The desired bound \eqref{contrapos} is now a direct consequence 
of the following elementary result: 

\noindent{\bf Claim:} {\it If $\mu \ge 1$, $\Psi \ge 1$, and $C \ge 1$ 
satisfy $\mu \le \Psi(C+\log\mu)$, then}
\begin{equation}\label{claim1}
  \mu \,\le\, \Psi(C + 2\log2 + 2\log\Psi + \log C)~.
\end{equation}
Indeed, the hypothesis implies that $\mu \le \mu_0$, where 
$\mu_0 \ge 1$ is uniquely determined by the relation $\mu_0 =
\Psi(C+\log\mu_0)$. Since $\log\mu_0 \le \sqrt{\mu_0}$, we have
$\mu_0 \le \Psi(C+\sqrt{\mu_0})$, hence
\[
  \sqrt{\mu_0} \,\le\, \frac12\Bigl(\Psi + \sqrt{\Psi^2 + 4C\Psi}
  \Bigr) \,\le\, \Psi + \sqrt{C\Psi}~.
\]
Using this bound and the fact that $\mu \le \Psi(C+\log\mu_0)$, we
easily obtain \eqref{claim1}. This concludes the proof of the 
claim, hence of Lemma~\ref{pseudo}. \QED

\medskip We now briefly recall the comparison argument with
semiclassical models which yields the lower bound \eqref{hatHbdd}.

\medskip\noindent
{\bf Proof of Lemma~\ref{hatH}.} Let $V(x,\epsilon) = x^2 + 
\frac{f'(x)^2}{\epsilon^2}$. By the Min-Max principle, it is 
sufficent to find a constant $M_3 \ge 1$ such that
\begin{equation}\label{eq.minmax1}
  \langle\hat H_\epsilon u\,,\,u\rangle \,=\, \int_\R \Bigl(
  |\partial_x u|^2 + V(x,\epsilon)|u|^2 \Bigr)\d x
  \,\ge\, \frac{\|u\|^2}{M_3\,\epsilon^{2\nu}}~, \quad 
  \hbox{for all }u\in \DD~,
\end{equation}
and for any $\epsilon\in (0,1]$ a nonzero function $\varphi_\epsilon\in 
C_0^\infty(\R) \subset \DD$ such that
\begin{equation}\label{eq.minmax2}
  \langle \hat H_\epsilon\varphi_\epsilon\,,\,\varphi_\epsilon\rangle 
  \,\le\, \frac{M_3}{\epsilon^{2\nu}}\,\|\varphi_\epsilon\|^2~.
\end{equation}
Hypothesis~\ref{hypf} ensures that there exists $L > 0$ such that
\[
  f'(x)^{2} \,\ge\, \frac{k^2}{2|x|^{2(k+1)}} \quad
  \hbox{for}\quad |x| \,\ge\, L~,
\]
and that $f$ has a finite set $\{x_1,\dots,x_N\}$ of (nondegenerate)
critical points. Hence there exists a partition of unity
$\sum_{j=0}^{N}\chi_j^2 = 1$ such that $\chi_j\in C_0^\infty(\R)$ for
$j \in \{1,\dots,N\}$, $\supp \chi_0\subset (-\infty,-L)
\cup(L,+\infty)$, and such that $f$ has exactly one critical point
(namely, $x_j$) in the support of $\chi_j$ for $j \in \{1,\dots,N\}$. 
The IMS localization formula of \cite{CFKS} provides a constant $C >
0$ such that
\begin{equation}\label{eq.IMSvr}
  \int_\R \Bigl(|\partial_x u|^2 + V(x,\epsilon)|u|^2\Bigr)\d x
  \,\ge\, \sum_{j=0}^N \int_\R \Bigl(|\partial_x u_j|^2 +
  V(x,\epsilon)|u_j|^2\Bigr)\d x -C \|u\|^2~,
\end{equation}
with $u_j = \chi_j u$. By construction, for $j\in \{1,\dots,N\}$, 
there exists $c_j > 0$ such that
\[
   V(x,\epsilon) \,\ge\, \frac{f'(x)^2}{\epsilon^2} \,\ge\,  
   \frac{c_j ^2(x-x_j)^2}{\epsilon^2}~, \quad \hbox{for all }
   x\in \supp \chi_j~.
\]
Hence the comparison with the harmonic oscillator hamiltonian
$-\partial_x^2 + c_j^2\,\epsilon^{-2}(x-x_j)^2$ implies that 
$\langle\hat H_\epsilon u_j\,,u_j\rangle \ge c_j\,\epsilon^{-1}\|u_j\|^2$
for all $j\in \{1,\dots,N\}$. For $j=0$, we notice that
\[
  V(x,\epsilon) \,\ge\, x^2+\frac{k^2}{2\epsilon^2|x|^{2(k+1)}} 
  \,\ge\, c_0\,\epsilon^{-\frac{2}{k+2}}~, \quad  
  \hbox{for all }x\in \supp \chi_0 \subset \R\setminus[-L,L]~.
\]
Thus the first term $\langle \hat H_\epsilon u_{0}\,,u_{0}\rangle$
is bounded from below by $c_{0}\,\epsilon^{-\nu}$. Summing up all the 
lower bounds of the terms in \eqref{eq.IMSvr}, and recalling 
that $\|u\|^2 \le \langle\hat H_\epsilon u_j\,,u_j\rangle$, 
we obtain \eqref{eq.minmax1} provided $M_3^{-1}(1+C) \le \min\{c_j\,;
\,0\le j\le N\}$. For the upper bound \eqref{eq.minmax2}, it 
suffices to take $\varphi_\epsilon(x) = \varphi(\epsilon^{\nu}x)$
with $\varphi \in C_0^\infty((1,2))$ and $\varphi \not\equiv 0$.
\QED

\medskip
Following exactly the same lines, one can also prove the following 
result, which is needed in Section~\ref{optimal}: 

\begin{lemma}\label{betamodif}
Let~$\beta : \R \to \R_+$ be the function defined in 
Section~\ref{optimal} and depicted in Fig.~1. There is a 
constant~$M_0 > 0$ such that, for any $u \in \DD$ and any 
$\epsilon \in (0,1]$, 
\begin{equation}\label{lowerbound2}
  \int_\R \Bigl(|\partial_x u|^2 + x^2 |u|^2 + \frac{\beta}{\epsilon}
  \,f'(x)^2 |u|^2\Bigr)\d x \,\ge\, \frac{M_0}{\epsilon^{\bar \nu}} 
  \|u\|_{L^2}^2~, \quad \hbox{where}\quad \bar \nu \,=\, 
  \frac2{k+4} ~\cdotp
\end{equation}
\end{lemma}

\proof
We assume that the parameter $A > 0$ entering the definition 
\eqref{betadef2} of the function $\beta$ is large enough so that
all critical points of $f$ are contained in the interval $[-A+1,A-1]$, 
and so that
\begin{equation}\label{Acond1}
  f'(x)^2 \,\ge\, \frac{k^2}{2|x|^{2(k+1)}}~, \quad \hbox{whenever}~ 
  |x| \ge A~.
\end{equation}
We now consider the potential $V(x,\epsilon) = x^2+\frac{\beta(x)f'(x)^2}
{\epsilon}$ and we introduce like in the proof of Lemma~\ref{hatH} 
a partition of unity $\sum_{j=0}^N\chi_j^2 \equiv 1$ with $\chi_j\in
C_0^\infty(\R)$ for all $j\in \{1,\dots,N\}$, $\supp \chi_0\subset
(-\infty,-A)\cup(A,+\infty)$, and such that $f$ has exactly one
critical point in $\supp \chi_j$ for $j\in \{1,\dots,N\}$. Again the 
IMS localization formula provides a constant $C > 0$ such that 
\[
  (1+C)\int_\R \Bigl(|\partial_x u|^2 + V(x,\epsilon)|u|^2\Bigr)
  \d x \,\ge\, \sum_{j=0}^N \int_\R \Bigl(|\partial_x u_j|^2 + 
  V(x,\epsilon) |u_j|^2\Bigr)\d x~,
\]
with $u_j = \chi_j u$. For $j\in \{1,\dots, N\}$ we find as in 
the proof of Lemma~\ref{hatH} a constant $c_j > 0$ such that
\begin{equation}\label{lastone}
  \int_\R \Bigl(|\partial_x u_j|^2 + V(x,\epsilon)|u_j|^2\Bigr)
  \d x \,\ge\, \frac{c_{j}}{\epsilon^{1/2}}\,\|u_j\|^2~.
\end{equation}
For $j=0$, we use the definition \eqref{betadef2} and the condition
\eqref{Acond1} to find an appropriate lower bound on $V(x,\epsilon)$. 
If $c_0 > 0$ is sufficiently small, then for any $x\in \supp \chi_0$
we have either
\[
  A \le |x|\le B_\epsilon~, \quad \hbox{so that}\quad
  V(x,\epsilon) \,\ge\, x^2 + \frac{\beta_0 k^2}{2 A^{2k}\epsilon 
  x^2} \,\ge\, c_0\,\epsilon^{-1/2}~,
\]
or 
\[ 
  |x| \ge B_\epsilon~, \quad \hbox{so that}\quad
  V(x,\epsilon) \,\ge\, x^2 + \frac{\beta_0\,\epsilon^{
  -\frac{k}{k+4}}k^2}{2\epsilon |x|^{2(k+1)}} \,\ge\, c_0
  \,\epsilon^{-\frac{2}{k+4}}~.
\]
Since $\bar \nu = \frac{2}{k+4} < \frac12$, it follows that 
$V(x,\epsilon) \ge c_0\,\epsilon^{-\bar\nu}$ if $x \in \supp\chi_0$, 
hence \eqref{lastone} holds for $j = 0$ too if we replace the exponent
$1/2$ by $\bar \nu$. We conclude as before. \QED


\end{document}